\documentclass[11pt, reqno]{article}
\usepackage{amsmath}
\usepackage{amsfonts}
\usepackage{graphicx}
\newtheorem{teo}{Theorem}[section]

\newtheorem{prop}[teo]{Proposition}
\newtheorem{cor}[teo]{Corollary}
\newtheorem{lemma}[teo]{Lemma}
\newtheorem{rem}[teo]{Remark}
\numberwithin{equation}{section}
\def\proof{{\bf Proof:}\ }
\def\endproof{\hfill $\Box$\par\vskip3mm}

\def\eq#1{(\ref{#1})}
\def\pr#1{\ref{#1}}
\def\neweq#1{\begin{equation}\label{#1}}
\def\endeq{\end{equation}}

\def\eps{\varepsilon}
\def\phi{\varphi}

\def\di{\displaystyle}
\def\ri{\rightarrow}

\def\ria{\rightharpoonup}

\def\rr{\mathbb{R}}
\def\nn{\mathbb{N}}
\def\mb{{\mathcal M}_b(\Omega)}

\def\demi{\frac{1}{2}}
\def\l1{L^1(\Omega)}
\def\lq{L^q(\Omega)}

\def\w1{W^{1,1}(\Omega)}

\def\crit{\frac{N+2}{N+1}}

\def\ligne{\noindent \rule{\textwidth}{0.3mm}}

\begin{document}
\title{\bf  Large time behaviour for a viscous
  Hamilton-Jacobi equation with Neumann boundary condition}
\author{\Large Sa\"{\i}d Benachour\footnote{Corresponding
    author. E-mail address: Said.Benachour@iecn.u-nancy.fr, Fax: 0033/383684534},  and Simona
  Dabuleanu\footnote{E-mail address: Simona.Dabuleanu@iecn.u-nancy.fr
}}
\date{}
\maketitle
\pagestyle{plain}

\vspace{-1cm}
\begin{center}{\small\em  Institut Elie Cartan  UMR 7502  UHP-CNRS-INRIA\\
BP 239  F-54506 Vandoeuvre-l\`es-Nancy  France}
\end{center}

\ligne

\noindent{\bf Abstract}\\

{\small We prove the existence and the uniqueness of strong solutions for the
viscous Hamilton-Jacobi equation:
  $u_t-\Delta u=a|\nabla u|^p,~t>0,~x\in\Omega~$ with Neumann boundary
  condition, and initial data
$\mu_0$, a continuous function. The domain $\Omega$ is a bounded
and convex open set with smooth boundary, $a\in\rr, a\neq 0$ and
$p>0$. Then, we study the large time behavior of the solution and
we show that for $p\in(0,1)$, the extinction in finite time of the
gradient of the solution  occurs, while for $p\geq 1$ the solution
converges uniformly to a
constant, as $t\to\infty$. \\

}

\ligne

\pagestyle{plain}

\section{Introduction and main results}
Consider the following initial boundary value problem:
\begin{equation}\label{part3.1.1}
\left\{\begin{array}{lll}
\di\frac{\partial u}{\partial t}-\Delta u=a|\nabla u|^p &\mbox{in}& (0,+\infty)\times\Omega,\\
\di\frac{\partial u}{\partial\nu}(t,x)=0&\mbox{on}&(0,+\infty)\times\partial\Omega,\\
u(0,x)=\mu_0(x)&\mbox{in}&\Omega,
\end{array}\right.
\end{equation}
where $a\in\rr$, $a\neq 0$, $p>0$ and $\Omega\subset\rr^N$ is a bounded
open set with smooth boundary of $C^3$ class.

The Cauchy problem in the whole space $\rr ^N$ has been
intensively studied (see \cite{Am,B,BSW,BL, BRV, GGK, SW}). As
well, in bounded domains $\Omega\subset\rr^N $, existence and
uniqueness results of the solutions for the Cauchy-Dirichlet
problem have been obtained in \cite{Al, BD, CLS,S}. In particular
the large time behavior of the solution to the Cauchy problem has
been analysed in \cite{BK,BLS,BLSS}, as $a<0$ and for initial data
$\mu_0$ a positive function. Thus, in \cite{BLSS}, we can find the
following result: if $a<0$, $p\in(0,1)$ and the initial data
$\mu_0$ is a periodic function, the extinction in finite time of
the solution of problem \eq{part3.1.1} occurs. Since, any positive
solution of the Cauchy problem is a super-solution of the
homogeneous Cauchy-Dirichlet problem, the result of  \cite{BLSS},
remain valid also in bounded domains for the Cauchy-Dirichlet
problem.

With respect to the Cauchy-Neumann we mention the results given in
\cite{D}, regarding the existence, uniqueness and regularity of
weak solutions, for $p\in(0,2)$, $a\in\rr, a\neq 0$ and initial
data $\mu_0$ a bounded Radon measure or a measurable function in
$\lq, \ q\geq 1$. To our knowledge the problem \eq{part3.1.1} has
not been investigated for the super-quadratic case, $p\geq 2$.

In this paper we consider the problem \eq{part3.1.1} when $\Omega$
is a bounded and convex open set, and we give some existence and
uniqueness results of the solutions when the initial data is a
continuous function in $\overline{\Omega}$. Then we study the
large time behavior of the solutions according to the exponent
$p$. The results rely on some remarkable estimates for the
gradient of the solutions of problem \eq{part3.1.1}, obtained by
using a Bernstein technique. These estimates, given in Theorem
\pr{part3.t.1} are used as the key argument in the proof of the
extinction result in Theorem
 \pr{part3.t.2}. More exactly we show that: if
$p\in(0,1)$ then, for any solution $u$ of problem \eq{part3.1.1} with initial
data in  $C(\overline{\Omega})$ there exists
$T^{*}>0$ and $c\in\rr$ such that:
$$u(t,x)\equiv c,\ \mbox{for all}\ t>T^{*}\ \mbox{and}\ x\in\Omega.$$
This property is called: ``the extinction of the gradient of the
solution $u$ in finite time``. Also, in Theorem  \pr{part3.t.2} we
prove that, for $p\geq 1$ any solution of problem \eq{part3.1.1}
converges uniformly to a constant,
as $t\to\infty$.\\

The notations used are mostly standard for the parabolic equations theory:\\
For   all $0<\tau<T\leq\infty$ we denote by $Q_T=(0,T)\times\Omega~ $,
$\Gamma_T=(0,T)\times\partial\Omega~$, $Q_{\tau,T}=(\tau,T)\times\Omega~$
and  $\Gamma_{\tau,T}=(\tau,T)\times\partial\Omega$.
$C(\overline{\Omega})$ is the space of continuous functions  on
$\overline{\Omega}$.
$C_b(\Omega)$ is the space of bounded continuous functions on  $\Omega$.
 $C_0(\overline{\Omega})$ the space of continuous functions on
 $\overline{\Omega}$ which vanish on the boundary $\partial\Omega$.
$C_c^{\infty}(\Omega)$ (
resp. $C_c^{\infty}(Q_T)$) the space of infinitely differentiable functions
on  $\Omega$ (resp. $Q_T$) with compact support in $\Omega$
(resp. $Q_T$).
$C^{0,1}([0,T)\times\overline{\Omega})$ is the space of continuous functions  $u$ on
$[0,T)\times\overline{\Omega}$  which are differentiable with respect to
$x\in\Omega$ and the derivatives $(\frac{\partial
  u}{\partial x_i})_{1\leq i\leq N} $  are in $C([0,T)\times\overline{\Omega})$.
$C^{1,2}(Q_T)$ is the space of continuous functions $u$ on $Q_T$ such that
the derivatives: $ \frac{\partial u}{\partial t},$
$(\frac{\partial u}{\partial x_{i}})_{1\leq i\leq N},$ $(\frac{\partial
  u}{\partial x_i\partial x_j})_{1\leq i,j\leq N },$  exist and belong to
$C(Q_T)$.
Suppose that  $\alpha$ is a positive real number and  $[\alpha]$ the
integer part of $\alpha$ such that $[\alpha]<\alpha$ then:
$C^{\alpha}(\overline{\Omega})$  and $C^{\alpha/2,\alpha}(\overline{Q})$
denote the usual H\"older spaces on the bounded open sets $\Omega\subset\rr^N$ and $Q\subset\rr^{N+1}$
respectively (for the  definitions see \cite{F,LSU}).\\
We denote by
$\mb$ the space of bounded Radon measures on
$\Omega$ endowed with the usual norm $\|~\|_{\mb}$.
For  $q\geq 1$,  $\|~ \|_q$ is the usual norm of the Lebesgue space
$L^q(\Omega)$.
$W^{1,q}(\Omega)$, $W^{1,q}(Q_T)$ and $W^{1,2}_q(Q_T)$ are the usual
Sobolev spaces in $\Omega$ respectively $Q_T$ (for the definitions see \cite{L}).\\
We denote by $(S(t))_{t\geq 0}$ the semigroup of contraction in
$L^q(\Omega),q\geq 1$ related to the heat equation with homogeneous Neumann
boundary condition (see \cite{Ro}). As we can see in \cite{D} this semigroup can be
extended, in a natural way, to the space of bounded Radon measures, $\mb$.

First we recall an existence and uniqueness result for the
solutions of problem \eq{part3.1.1} when $p\in(0,2)$ (for further
details see \cite{D}).
\begin{teo}\label{part3.t.5}\cite{D}
Let $\mu_0\in\mb$. Then \eq{part3.1.1} admits a weak solution
$u\in L^{\infty}(0,T;\l1)\cap L^1(0,T;\w1)\cap
C^{1+\delta/2,2+\delta}(\overline{Q_{\tau,T}})$,
$T>0,\tau\in(0,T)$, $\delta\in(0,1)$, such that $|\nabla u|^p\in
L^1(Q_T)$, in the following cases:
\begin{itemize}
\item[$i)$]$0<p<2/(N+1)$. The solution is unique if $\Omega$ is
convex.

\item[ $ ii)$] $2/(N+1)\leq p<1$. The solution is unique if
$\mu_0\in L^q(\Omega)$ for some $q>pN/(2-p)$ and $\Omega$ a convex
open set.

\item[$iii)$] $1\leq p<(N+2)/(N+1)$. The solution is unique. If
$\mu_0\in L^q(\Omega)$ for some $q\geq 1$ then $u\in
C([0,T];\lq)\cap L^p(0,T; W^{1,pq}(\Omega))$.

\item[$iv)$] $(N+2)/(N+1)\leq p<2$, $\mu_0\in \lq$ and
$q>q_c=\frac{N(p-1)}{2-p}$. There holds $u\in C([0,T];\lq)\cap
L^p(0,T; W^{1,pq}(\Omega))$ and the solution is unique in this
space.
\item[$v)$]$(N+2)/(N+1)\leq p<2$, $\mu_0\in L^1(\Omega)$,
$\mu_0\geq 0$.
\end{itemize}
Moreover, this solution satisfies \eq{part3.1.1} in the mild
sense:
$$u(t)=S(t)\mu_0+a\int\limits_0^tS(t-s)|\nabla u(s)|^p\,ds,\ \ t\in(0,T).$$
\end{teo}

 In Theorem \pr{part3.t.1} below we prove the existence and the  uniqueness of
 solutions of  problem \eq{part3.1.1}, for $p>0$, $\Omega$ a bounded and
 convex open set with smooth boundary, and for initial data
$\mu_0\in C(\overline{\Omega})$. We give also some gradient
estimates of the solution $u$ of problem \eq{part3.1.1} which will
be very useful in the proof of Theorem \pr{part3.t.2}.

Let $u$ be a function in $C(\overline{Q_{\infty}})$. For any
$t\geq 0$ denote by:
\begin{equation}\label{part3.1.2}
M(t)=\max\limits_{x\in\overline{\Omega}}u(t,x)
\end{equation}
and
\begin{equation}\label{part3.1.3}
m(t)=\min\limits_{x\in\overline{\Omega}}u(t,x),
\end{equation}
\begin{teo}\label{part3.t.1}
Consider $a\in\rr, a\neq 0$, $p>0$ and $\mu_0\in
C(\overline{\Omega})$, where $\Omega$ is a bounded and convex open
set. Then, the problem \eq{part3.1.1} admits a unique solution:
 $$u\in
C(\overline{Q_T})\cap C^{1+\delta/2,2+\delta}(\overline{Q_{\tau,T}})$$
for any  $T>0$ and $\tau\in(0,T)$. Moreover, we have:
\begin{equation}\label{part3.1.4}
t\ri M(t)\  \mbox{is a decreasing function in}\ \rr,
\end{equation}
\begin{equation}\label{part3.1.5}
t\ri m(t)\ ~\mbox{ is a non-decreasing function in}\ \rr ,
\end{equation}
\begin{equation}\label{part3.1.5.1}
\|\nabla
u(t)\|_{\infty}\leq\left(\demi\right)^{1/2}(M(s)-m(s))(t-s)^{-\demi}\
\mbox{for all}\  t>s\geq 0,
\end{equation}
 and for $p\neq 1$
\begin{equation}\label{part3.1.6}
\|\nabla u(t)\|_{\infty}\leq
\left(\frac{\max\{p,2\}}{ap|1-p|}\right)^{1/p}(M(s)-m(s))^{1/p}(t-s)^{-1/p}\
\mbox{for all}\  t>s\geq 0.
\end{equation}
\end{teo}

For the proof we are using the Bernstein technique. This method
can be found in \cite{BL,CLS,GGK} and \cite{Li2},  where formulas
similar to \eq{part3.1.5.1} and \eq{part3.1.6} are obtained for
the Cauchy problem in $\rr^N$. This method has also been used by
Ph. Benilan \cite{Be} in order to obtain remarkable estimates for
the solutions of ``the porous medium equation''

In the next result we are going to analyze the large time behavior
of the solutions for problem \eq{part3.1.1}.
\begin{teo}\label{part3.t.2}
Consider $a\in\rr, a\neq 0$, $p>0$ and $\Omega$ a bounded and
convex domain. Let $\mu_0\in C(\overline{\Omega})$ and denote by
$u$ a solution of problem \eq{part3.1.1} corresponding to $\mu_0$.
Then:
\begin{itemize}
\item[$i)$] If $p\in(0,1)$, the extinction  of the gradient of
  $u$ in finite time occurs, in other words: \\
there exists $T^*\in[0,+\infty)$ and  $c\in\rr$ such that:
$$u(t,x)\equiv c\ \mbox{for all}\ t\geq T^{*} \ \mbox{and}\ x\in\overline{\Omega},$$
\item[$ii)$] If $p\in[1,+\infty)$, then $u(t,\cdot)$ converges
  uniformly on $\overline{\Omega}$ to a constant, as $t\to\infty$. Moreover
  the decreasing rate is given by:
$$
M(t)-m(t)\leq \left(\frac{8f(t/2)}{t^2}\right)^{\frac{1}{\gamma}},\ \
\forall\ t>0,$$
 where $f$ is defined in \eq{part3.5.101} and $\gamma$ in \eq{part3.4.2}.
\end{itemize}
\end{teo}
\begin{rem}
From \eq{part3.5.101} we have:
$$\left(\frac{8f(t/2)}{t^2}\right)^{\frac{1}{\gamma}}=\left\{\begin{array}{lll}
C_1t^{-\frac{2}{\gamma}}e^{-C_2t}&\mbox{if}&p=1,\\
C_3t^{-\frac{\alpha+1}{\gamma(\alpha-1)}}&\mbox{if}&p>1.\end{array}\right.$$
where $\gamma$ and $\alpha $ are given by \eq{part3.4.2} and \eq{part3.4.5},
and $C_1, C_2, C_3$ are positive constants which depends only on $p, N,
\gamma$.
\end{rem}

The proof of Theorem \pr{part3.t.2} follows the same ideas as in
\cite{BLSS}. In this paper, the authors investigate the large time
behaviour for the Cauchy problem in the whole space $\rr^N$ and
for initial data periodic functions. We mention that the key
arguments of the proof are the relations \eq{part3.1.5.1} and
\eq{part3.1.6} above.
\begin{rem}\label{part3.r.0}
Theorem \pr{part3.t.2} is valid for any $a\in\rr, a\neq 0$, while in
\cite{BLS} and \cite{BLSS} the result is proved for $a<0$.
\end{rem}
The next result is a simple consequence of Theorem \pr{part3.t.5} and
Theorem \pr{part3.t.2}  above.
\begin{cor}\label{part3.c.1}
Let  $\Omega$ be a bounded and convex domain with smooth boundary and let
$a\in\rr, a\neq 0$. Then:
\begin{itemize}
\item[$i)$]
If $p\in(0,1)$ and $\mu_0$ is a bounded Radon measure, the
extinction in finite time of the gradient of any weak solution $u$ of problem \eq{part3.1.1}  occurs.
\item[$ii)$] The weak solution $u(t,\cdot)$ of problem  \eq{part3.1.1}
  converges uniformly in  $\overline{\Omega}$, to a constant $c\in\rr$,
  as $t\to\infty$, in the two cases below:\\
$a)$ $p\in[1,\crit)$ and $\mu_0$ is a bounded Radon measure, \\
$b)$ $p\in[\crit,2)$  and $\mu_0\in \lq$,
  $q>q_c=\frac{N(p-1)}{2-p}$.
\end{itemize}
\end{cor}

This paper is organized as follows: In section $2$, we give some
preliminary results. In section $3$, we introduce the technique of
Bernstein to obtain some uniform estimates for the gradient of the
solution of problem \eq{part3.1.1} and we prove the Theorem
\pr{part3.t.1}.Finally section $4$ is devoted to the proof of
Theorem \pr{part3.t.2}, which concerns the large time behaviour of
solutions.

\section{Preliminary results}
We start with some auxiliary results.
\begin{lemma}\label{part3.l.1}
Let $\Omega\subset\rr^N$ be a bounded domain with smooth boundary and
consider $\mu_0\in C(\overline{\Omega})$. Denote by
$m=\min\limits_{x\in\Omega}\mu_0(x)$ and
$M=\max\limits_{x\in\Omega}\mu_0(x)$. \\
Then, there exists a sequence
$(u_0^n)_{n\geq 1}\subset C^{3+\beta}(\overline{\Omega}), (\beta\in(0,1))$
such that:
\begin{equation}\label{part3.2.1}
u_0^n\searrow \mu_0\ \mbox{as}\ n\to\infty,
\end{equation}
\begin{equation}\label{part3.2.1.1}
\ m+\frac{1}{2^{n+1}}\leq u_0^n\leq
M+\frac{1}{2^{n-1}},\ \forall n\geq 1,
\end{equation}
and
\begin{equation}\label{part3.2.2}
\frac{\partial u_0^n}{\partial \nu}=0\ \mbox{on}\  \partial\Omega.
\end{equation}
\end{lemma}
\proof\\
For any $n\in\nn^*$, denote by $v_0^n=\mu_0+\frac{1}{2^n}$, then
  $(v_0^n)_n\subset C(\overline{\Omega})$.
For $t>0$ let us set:
$$v^n(t)=S(t)v_0^n.$$
 Then $v^n\in
C(\overline{Q_{\infty}})\cap
C^{\infty}(\overline{Q_{\tau,\infty}})$ for all
$\tau\in(0,\infty)$, and:
$$\frac{\partial v^n}{\partial\nu}(t,x)=0\ \ \ \mbox{for all}\ (t,x)\in\Gamma_T.$$
Since $v^n\in C(\overline{Q_{\infty}})$, there exists $t_n$ close enough
from $0$ such that:
\begin{equation}\label{part3.2.4}
|v^n(t_n,x)-v_0^n(x)|<\frac{1}{2^{n+2}},\ \ \ \ \forall  x\in\Omega.
\end{equation}
Denote by: $$u_0^n(x)=v^n(t_n,x),\ \ \ \ x\in\Omega.$$ Then $u_0^n\in
C^{\infty}(\overline{\Omega})$ and satisfies condition
\eq{part3.2.2}. Moreover, thanks to \eq{part3.2.4} we have on the one hand:
$$u_0^n-\mu_0=(u_0^n-v_0^n)+(v_0^n-\mu_0)\leq\frac{1}{2^{n+2}}+\frac{1}{2^n}\leq
\frac{1}{2^{n-1}},$$
on the other hand:
$$u_0^n-\mu_0=(u_0^n-v_0^n)+(v_0^n-\mu_0)\geq -\frac{1}{2^{n+2}}+\frac{1}{2^{n}}\geq\frac{1}{2^{n+1}}.$$
which yields \eq{part3.2.1.1}.\\
To prove that $(u_0^n)_n$ is a decreasing sequence, let compute:
$$u_0^n-u_0^{n+1}=(u_0^n-v_0^n)+(v_0^n-v_0^{n+1})+(v_0^{n+1}-u_0^{n+1})\geq
-\frac{1}{2^{n+2}}+\frac{1}{2^{n+1}}-\frac{1}{2^{n+3}}=\frac{1}{2^{n+3}}>0.$$
And finally we obtain \eq{part3.2.1}.
\endproof
\begin{lemma}\label{part3.l.2}
Let $\Omega\subset\rr^N$ be a convex and bounded domain and $q$ a
real number such that $q>N$. From the Sobolev embedding, $
W^{1,q}(\Omega)\hookrightarrow C(\overline{\Omega})$, for all
$u\in W^{1,q}(\Omega)$, the following quantities:
$$M_u=\max\limits_{x\in\overline{\Omega} }u(x)\ \mbox{ and}
 \  m_u=\min_{x\in\overline{\Omega}}u(x),$$
are well defined. Moreover we have :
\begin{equation}\label{part3.2.6}
M_u-m_u\leq C\|\nabla u\|_q,
\end{equation}
where $C$ is a positive constant depending only on $q$, $N$ and $\Omega$.
\end{lemma}
\proof
 The proof is similar to that of Lemmas 7.16 and 7.17
in \cite{GT}.  $\Omega$ being a convex set, for all $x,y\in\Omega$
we have
 $(1-t)x+ty\in\Omega$ for any $t\in[0,1]$. Let $u\in W^{1,q}(\Omega)$, then:
$$u(x)-u(y)=\int\limits_0^1\nabla u((1-t)x+ty)\cdot (x-y)\, dt,$$
which yields:
$$u(x)-\frac{1}{|\Omega|}\int\limits_{\Omega}u(y)\,
dy=\frac{1}{|\Omega|}\int\limits_{\Omega}\int\limits_{0}^{1}\nabla
u((1-t)x+ty)\cdot (x-y)\,dt\,dy.$$
Denote by:
$$u_{\Omega}=\frac{1}{|\Omega|}\int\limits_{\Omega}u(y)\,
dy\ \ \ \ \
\mbox{and}\ \ \ \ \ d=diam(\Omega).$$
Then,
$$\begin{array}{l}
|u(x)-u_{\Omega}|\leq\di\frac{d}{|\Omega|}\int\limits_{\Omega}
\int\limits_{0}^1|\nabla u((1-t)x+ty)|dtdy
\leq\di\frac{d}{|\Omega|}\int\limits_0^1\int\limits_{\Omega}|\nabla
u((1-t)x+ty)|dydt.
\end{array}$$
We replace $(1-t)x+ty=\zeta$ and, for any
$t\in[0,1]$, we denote by $\Omega_t$ the set:
$$\Omega_t=\{\zeta=(1-t)x+ty; y\in\Omega\}\subset\Omega,$$
then:
$$|u(x)-u_{\Omega}|\leq\frac{d}{|\Omega|}\int\limits_0^1\int\limits_{\Omega_t}|\nabla
u(\zeta)|t^{-N} dt\,d\zeta.$$
Using the H\"older inequality for $q>N$ we get:
$$\begin{array}{l}
\di|u(x)-u_{\Omega}|\leq\frac{d}{|\Omega|}\int\limits_{0}^1\Big[\int\limits_{\Omega}|\nabla
u(\zeta)|^q\, d\zeta\Big]^{1/q}|\Omega_t|^{1-1/q}t^{-N}\,dt\\
\ \ \ \ \ \leq\di\frac{d}{|\Omega|}\|\nabla u\|_{q}\cdot\int\limits_0^1
t^{N(1-1/q)}t^{-N}|\Omega|^{1-1/q}\,dt\\
\ \ \ \ \ \leq\di\frac{d}{|\Omega|^{1/q}}\|\nabla u\|_{q}\cdot\int\limits_0^1
t^{-N/q}\,dt
\leq\di\frac{d}{|\Omega|^{1/q}}\frac{q}{q-N}\|\nabla u\|_q.
\end{array}$$
Finally, for $x,y\in\Omega$ we obtain:
$$|u(x)-u(y)|\leq|u(x)-u_{\Omega}|+|u_{\Omega}-u(y)|\leq\frac{2d}{|\Omega|^{1/q}}\cdot\frac{q}{q-N}\|\nabla
  u\|_q,$$
and relation \eq{part3.2.6} follows. Thus, the Lemma  \pr{part3.l.2} is achieved.
\endproof
\begin{lemma}\label{part3.l.3}
Let $\Omega\subset\rr^N$ be a convex and bounded domain, then for
$u\in C^2(\overline{\Omega})$ such that $\di\frac{\partial
  u}{\partial\nu}\mid_{\partial\Omega}=0$, we have:
$$\frac{\partial}{\partial \nu}|\nabla u|^2\leq 0\ \mbox{on}\ \partial\Omega.$$
\end{lemma}
For the proof see Lemma I.1, p. 350 in \cite{Li2}.\\

The following lemma is a comparison principle for parabolic
nonlinear equations, which generalize the result obtained in
\cite{Kp}, to less regular functions.
\begin{lemma}\label{part3.l.4}
Let $\Omega\subset\rr^N$ be  a convex and bounded open set with smooth
boundary and denote by ${\mathcal N}$, the nonlinear parabolic operator,
defined by:
$${\mathcal N}(u)=\frac{\partial u}{\partial t}-\Delta u-f(t,x,u,\nabla u)$$
where $f$ is a uniformly continuous function satisfying:\\
 for all $r>0$ there exists $L_r>0$  such that:
\begin{equation}\label{part3.2.4.1}
\begin{array}{l}
|f(t,x,y_1,v_1)-f(t,x,y_2,v_2)|\leq L_r(|y_1-y_2|+|v_1-v_2|),\ \\
\qquad\qquad\qquad\mbox{for all}\ (t,x)\in Q_T \
\mbox{and}\ \ y_1,y_2\in(-r,r),\ v_1, v_2\in B_r(0),\end{array}
\end{equation}
where:
$$B_r(0)=\{\xi\in\rr^N;\  |\xi|<r\}.$$
Let $u^1$ and  $u^2$ be two functions in
$C^{0,1}(\overline{Q_{T}})\cap C^{1,2}(Q_T),$ such that:
\begin{equation}\label{part3.2.4.2}
\left\{\begin{array}{l}
{\mathcal N}(u^1)(t,x)\leq 0\leq {\mathcal N}(u^2)(t,x)\ \mbox{for all}
\ (t,x)\in Q_T\\
\di\frac{\partial u^1}{\partial \nu}\leq\frac{\partial u^2}{\partial
  \nu}\ \mbox{on}\ \Gamma_T\\
u^1(0,x)\leq u^{2}(0,x)\ \mbox{for all}\ x\in\Omega\end{array}\right.
\end{equation}
Then $$u^1\leq u^2\ \mbox{on}\ Q_T.$$
\end{lemma}
We begin the proof by the following useful remark:
\begin{rem}\label{part3.r.2.1}
Let $\Omega$ be a convex open set  in $\rr^N$ with smooth boundary $\partial\Omega$, which contains the origin.
For $x\in\partial\Omega$, denote by $\nu(x)$  the unit outward normal on
$\partial\Omega$ at the point $x$. Then:
$$x\cdot\nu(x)>0.$$
\end{rem}
{\bf Proof of Lemma \pr{part3.l.4}:}
Supposing first that $\Omega $ is a convex open set which contain the
origin and denoting by
$$R=\max\{\sup\limits_{(t,x)\in Q_T}|u_1(t,x)|;\sup\limits_{(t,x)\in
    Q_T}|\nabla u_1|(t,x);\sup\limits_{(t,x)\in Q_T}|u_2(t,x)|;\sup\limits_{(t,x)\in
    Q_T}|\nabla u_2|(t,x)\},$$
then, from \eq{part3.2.4.1}, there exists $L_R>0$  such that:
\begin{equation}\label{part3.2.100}
|f(t,x,u_1,\nabla u_1)-f(t,x,u_2,\nabla u_2)|\leq L_R(|u_1-u_2|+|\nabla
u_1-\nabla u_2|), \ \ \ \ (t,x)\in Q_T.
\end{equation}
For any $\eps\in(0,1)$ consider the function:
$$z(t,x)=u_1(t,x)-u_2(t,x)-\eps e^{Ct}(1+|x|^2)^{\demi},$$
where  $C=2L_R+N$. Then, using the regularity of $u_1$ and  $u_2$ we deduce
that:
 $$z\in C^{0,1}(\overline{Q_{T}})\cap C^{1,2}(Q_T).$$
For any $t\in[0,T]$ let us define the function:
$$\phi(t)=\max\{\sup\limits_{x\in\Omega}z(t,x);0\},$$
then $\varphi\in C([0,T])$, and for any $t\in(0,T]$ we can define:
$$\overline{\varphi'}(t)=\limsup\limits_{h\searrow
  0}\frac{\varphi(t)-\varphi(t-h)}{h}.$$
Thus, in order to prove that:
\begin{equation}\label{part3.2.4.3}
z\leq 0 \ \mbox{in}\ Q_T,
\end{equation}
 we need to show that:
\begin{equation}\label{part3.2.4.4}
\overline{\varphi'}(t)\leq L_R \varphi(t)\ \mbox{for all}\ \ t\in(0,T).
\end{equation}
Indeed, as $\varphi(0)=0$ and $\varphi\geq 0$, we can apply Theorem 4.1 in
\cite{H} to the differential inequality  \eq{part3.2.4.4}
and we deduce that
$\varphi\equiv 0$. Which implies \eq{part3.2.4.3}.\\
Proof of  \eq{part3.2.4.4}:
Consider $t\in(0,T]$.\\
There are two possibilities. Either $\varphi(t)=0$ and
\eq{part3.2.4.4} holds because, in this case,
$\overline{\varphi'}(t)\leq 0$. Or $\varphi(t)>0$, and in
particular, there exists $x_0\in\overline{\Omega}$ such that:
$$z(t,x_0)=\varphi(t)>0.$$
We claim that $x_0\notin\partial\Omega$. Indeed, if $x_0\in\partial\Omega$,
on the one hand:
$$\frac{\partial z}{\partial\nu}(t,x_0)=\lim\limits_{\lambda\nearrow
    0}\frac{z(t,x_0+\lambda\nu)-z(t,x_0)}{\lambda}\geq 0.$$
On the other hand, thanks to hypothesis \eq{part3.2.4.2} and to Remark
\pr{part3.r.2.1} we have:
$$\frac{\partial
  z}{\partial\nu}=\frac{\partial u_1}{\partial\nu}-\frac{\partial
  u_2}{\partial\nu}-\eps
e^{Ct}\frac{x\cdot\nu}{(1+|x|^2)^{\demi}}\leq -\eps
e^{Ct}\frac{x\cdot\nu}{(1+|x|^2)^{\demi}}<0\ \mbox{on}\ \Gamma_T.$$
So, we have a contradiction. Consequently, $x_0\in\Omega$ is a positive
maximum point for the function $\Omega\ni x\mapsto z(t,x)$. In particular
we have:
\begin{equation}\label{part3.2.4.5}
\nabla z(t,x_0)=0\ \ \mbox{and}\ \ \Delta z(t,x_0)\leq 0.
\end{equation}
Since, for any $h>0$, $z(t-h,x_0)\leq \varphi(t-h)$, we deduce:
\begin{equation}\label{part3.2.4.6}
\overline{\varphi'}(t)\leq\lim\limits_{h\nearrow
  0}\frac{z(t,x_0)-z(t-h,x_0)}{h}=\frac{\partial z}{\partial t}(t,x_0).
\end{equation}
On the other hand, thanks to \eq{part3.2.100} and
\eq{part3.2.4.5}, at $(t,x_0)$, we have:
\begin{equation}\label{part3.2.4.7}
\begin{array}{l}
\di\frac{\partial z}{\partial t}(t,x_0)=\frac{\partial u_1}{\partial
  t}(t,x_0)-\frac{\partial u_2}{\partial t}(t,x_0)-\eps
Ce^{Ct}(1+|x_0|^2)^{\demi}\\
~~~\\
 \di\leq \Delta (u_1-u_2)(t,x_0)+f(t,x_0,u_1,\nabla
u_1)-f(t,x_0,u_2,\nabla u_2) -\eps
Ce^{Ct}(1+|x_0|^2)^{\demi}\\
~~~\\
 \leq\Delta z+\eps
e^{Ct}\frac{(N+(N-1)|x_0|^2)}{(1+|x_0|^2)^{3/2}}+L_R|u_1-u_2|+L_R|\nabla
u_1-\nabla u_2| -\eps
Ce^{Ct}(1+|x_0|^2)^{\demi}\\
~~~\\
 \di\leq\eps
e^{Ct}\Big(N+L_R\frac{|x_0|}{(1+|x_0|^2)^{\demi}}\Big)+L_R
z(t,x_0)- \eps
(L_R+N)e^{Ct}(1+|x_0|^2)^{\demi}\\
~~~\\
 \di\leq L_R z(t,x_0)+\eps (N+L_R) e^{Ct}-\eps
(N+L_R) e^{Ct}(1+|x_0|^2)^{\demi}\leq L_R \varphi(t).
\end{array}
\end{equation}
Recall  that $C=2L_R+N$ and $z(t,x_0)=\varphi(t)$.
Combining \eq{part3.2.4.6} and \eq{part3.2.4.7} we deduce \eq{part3.2.4.4}. Thus
\eq{part3.2.4.3} holds. We may let $\eps\searrow 0$ in
\eq{part3.2.4.3} and we get: $$u_1\leq u_2\ \mbox{in}\ \ Q_T.$$

For the general case when $\Omega$ do not contains the origin, it
is possible to translate the problem on a domain which contains
the origin since the first equation of \eq{part3.1.1} is invariant
to the translation. For example we can carry the study of the
problem on $\ \Omega_{x_0}=\Omega-x_0$, where $x_0\in\Omega$.
\endproof

In the sequel we denote by $G:(0,+\infty)\times\Omega\times\Omega$ the heat
kernel for the homogeneous Neumann boundary value problem, then, for fix $y\in\Omega$, $G(\cdot,\cdot,y)$ verifies:
$$\left\{\begin{array}{l}
\di\frac{\partial G}{\partial t}(t,x,y)=\Delta_{x}G(t,x,y)\ \ \mbox{in}\
\ Q_{\infty},\\
\di\frac{\partial G}{\partial\nu}(t,x,y)=0\ \ \ ~\mbox{on}\ \ \Gamma_{\infty},\\
G(t,x,y)\underset{t\to 0}{\ria}\delta_y(x)\ \ \mbox{weakly in}\ \mb.
\end{array}\right.$$
The proof of the following property on the heat kernel can be found in
\cite{EI,F}.
\begin{lemma}\label{part3.l.5}\cite{EI, F}
Let $\Omega$ be a bounded open set with smooth boundary and $G$ the heat
kernel for the homogeneous Neumann boundary value problem.
Then for any $l\in\nn$ and $\alpha\in\nn^N$, and for any $T>0$, there
exists two positive constants $c>0$ and $C(T)>0$ such that:
\begin{equation}\label{part3.2.7}
|D_x^{\alpha} D_t^lG(t,x,y)|\leq C(T)
t^{-(\frac{N}{2}+\frac{|\alpha|}{2}+l)}e^{-c\frac{|x-y|^2}{t}}
\end{equation}
 for all $(t,x,y)\in(0,T)\times\Omega\times\Omega$.
\end{lemma}
Consider, $\mu_0\in L^{\infty}(\Omega)$ and $S(t)\mu_0$  the solution of
the heat equation with initial data $\mu_0$ and with  homogeneous
Neumann boundary condition. Then:
$$S(t)\mu_0(x)=\int\limits_{\Omega}G(t,x,y)\mu_0(y)\, dy.$$
Thanks to \eq{part3.2.7}, for any $l\in\nn$ and $\alpha\in \nn^N$ and for any
$T>0$ we have:
\begin{equation}\label{part3.2.8}
\|D_x^{\alpha} D_t^l S(t)\mu_0\|_{\infty}\leq C(T)\|\mu_0\|_{\infty}
t^{-(\frac{|\alpha|}{2}+l)},
\end{equation}
where $C(T)$ is a positive constant.

\section{Proof of Theorem \pr{part3.t.1}}
We prove the theorem for $a>0$. If $a<0$, then $-a>0$ and we notice that if
$v$ is the solution of problem \eq{part3.1.1} with initial data $-\mu_0$
instead of $\mu_0$ and $-a$ instead of
$a$, then $u=-v$ is the solution of problem \eq{part3.1.1} corresponding to
data $a$ and $\mu_0$.

The proof follows five steps:\\

{\bf First step: ``Smoothing''.}\\
Consider  $\mu_0\in C(\overline{\Omega})$ and denote by
$M(0)=\max\limits_{x\in\overline{\Omega}}\mu_0(x)$ and
$m(0)=\min\limits_{x\in\overline{\Omega}}\mu_0(x)$. Then, from Lemma
\pr{part3.l.1}, there exists a sequence of functions $(u_0^n)_{n\geq 1}$ satisfying:
\begin{equation}\label{part3.3.1}
\left\{\begin{array}{l}
u_0^n\searrow \mu_0\ \ \mbox{as}\ \ n\to\infty ,\\
\di\frac{\partial u_0^n}{\partial\nu}=0\\
m(0)+\frac{1}{n}\leq u_0^n\leq M(0)+\frac{2}{n},\ \ \forall n\geq 1.
\end{array}\right.
\end{equation}
As in \cite{BL} and \cite{GGK} we need to introduce a smooth
function, related to $\xi\to a|\xi|^p$.  So, for any
$\eps\in(0,1)$ we consider the application $F_{\eps}:\rr^N\to\rr$
defined by:
\begin{equation}\label{part3.3.2}
F_{\eps}(\xi)=\left\{\begin{array}{lll}
a(\eps+|\xi|^2)^{p/2}&\mbox{if}&0<p\leq 1,\\
a(-\eps+|\xi|^2)(\eps+|\xi|^2)^{\frac{p-2}{2}}&\mbox{if}&1<p<2,\\
a|\xi|^p&\mbox{if}&p\geq 2.\end{array}\right.
\end{equation}
With $\rho>0$ fixed, let us show that for any $\xi_1,\xi_2\in B_{\rho}(0)$
and $\eps\in(0,1)$ we have:
\begin{equation}\label{part3.3.3}
|F_{\eps}(\xi_1)-F_{\eps}(\xi_2)|\leq K \rho^{\max\{p-1,0\}}|\xi_1-\xi_2|^{\min\{p,1\}},
\end{equation}
where $K$ is a positive constant depending only on $p$ and $a$.

To prove  \eq{part3.3.3} we can distinguish among the three cases.
So, using the Mean Value Theorem there exists $\lambda\in[0,1]$ such that:\\
The case $0<p\leq 1$:\\
$$\begin{array}{l}
\di|F_{\eps}(\xi_1)-F_{\eps}(\xi_2)|=a[
(\eps+|\xi_1|^2)^{p/2}-(\eps+|\xi_2|^2)^{p/2}]\leq
a[(\eps+|\xi_1|^2)^{1/2}-(\eps+|\xi_2|^2)^{1/2}]^p\\
\qquad\di\leq
a\left(\frac{|\lambda\xi_1+(1-\lambda)\xi_2|}{(\eps+|\lambda\xi_1+(1-\lambda)\xi_2|^2)^{1/2}}|\xi_1-\xi_2|\right)^p\leq
a|\xi_1-\xi_2|^p.\end{array}$$
The case $1<p<2$:\\
$$\begin{array}{l}
\di|F_{\eps}(\xi_1)-F_{\eps}(\xi_2)|\leq |\nabla
F_{\eps}(\lambda\xi_1+(1-\lambda)\xi_2)\cdot(\xi_1-\xi_2)|\\
\qquad\di \leq a
\frac{2|\lambda\xi_1+(1-\lambda)\xi_2|(2\eps+\frac{p}{2}(|\lambda\xi_1+(1-\lambda)\xi_2|^2-\eps))}{(\eps+|\lambda\xi_1+(1-\lambda)\xi_2|^2)^{2-p/2}}|\xi_1-\xi_2|\\
\qquad\di\leq 4a|\lambda\xi_1+(1-\lambda)\xi_2|^{p-1}|\xi_1-\xi_2|\leq
4 a\rho^{p-1}|\xi_1-\xi_2|.\end{array}$$
The case $p\geq 2$:
$$\begin{array}{l}
|F_{\eps}(\xi_1)-F_{\eps}(\xi_2)|\leq |\nabla
F_{\eps}(\lambda\xi_1+(1-\lambda)\xi_2)||\xi_1-\xi_2|\\
\qquad\leq ap\cdot |\lambda\xi_1+(1-\lambda)\xi_2|^{p-1}|\xi_1-\xi_2|\leq ap\rho^{p-1}|\xi_1-\xi_2|.\end{array}$$

Moreover, $F_{\eps}\in C^{\infty}(\rr^{N})$ and satisfies the following inequalities:
\begin{eqnarray}
(\nabla F_{\eps})(\xi)\cdot\xi-F_{\eps}(\xi)\leq a(p-1)|\xi|^p \ \ \mbox{if}\
0<p\leq 1\label{part3.3.4},\\
(\nabla F_{\eps})(\xi)\cdot\xi-F_{\eps}(\xi)\geq a(p-1)|\xi|^p \ \ \mbox{if}\
\
p>1.\ \ \label{part3.3.5}
\end{eqnarray}
Indeed, when $0<p\leq 1$ we have:
$$\begin{array}{l}
(\nabla
F_{\eps})(\xi)\cdot\xi-F_{\eps}(\xi)=a\di\frac{p|\xi|^2-(\eps+|\xi|^2)}{(\eps+|\xi|^2)^{1-p/2}}=a\di\frac{(p-1)|\xi|^2-\eps}{(\eps+|\xi|^2)^{1-p/2}}\\
=-a\di\frac{\eps+(1-p)|\xi|^2}{(\eps+|\xi|^2)^{1-p/2}}=-a\di\frac{\eps+(1-p)|\xi|^2}{(\eps+|\xi|^2)}\cdot(\eps+|\xi|^2)^{p/2}\\
\leq -a\di\frac{(1-p)(\eps+|\xi|^2)}{\eps+|\xi|^2}\cdot(\eps+|\xi|^2)^{p/2}\leq
  -a(1-p)(\eps+|\xi|^2)^{p/2}\leq a(p-1)|\xi|^p.
\end{array}$$
If $1<p<2$ then:
$$\begin{array}{l}
\di(\nabla
F_{\eps})(\xi)\cdot\xi-F_{\eps}(\xi)=a\frac{(p-1)(\eps+|\xi|^2)^2+3\eps(2-p)|\xi|^2+\eps^2(2-p)}{(\eps+|\xi|^2)^{2-p/2}}\\
\qquad\qquad\qquad\di\geq a(p-1)(\eps+|\xi|^2)^{p/2}\geq a(p-1)|\xi|^p,
\end{array}$$
and finally, for $p\geq 2$ we have:
$$\di(\nabla
F_{\eps})(\xi)\cdot\xi-F_{\eps}(\xi)=ap|\xi|^{p-2}\xi\cdot\xi-a|\xi|^p=a(p-1)|\xi|^p.$$

For any $n\in\nn$, let denote by:
\begin{equation}\label{part3.3.5.1}
\rho_n=\sup\limits_{x\in\overline{\Omega}}\{|\nabla u_0^n|(x)\}.
\end{equation}
Then, there exists $\delta\in(0,1)$ and a function $F_{n,\eps}$,  such that:
\begin{eqnarray}
F_{n,\eps}\in C^{2+\delta}(\rr^N),\ \ \ \ \ \ \ \qquad\qquad \ \ \ \ \ \ \label{part3.3.5.5}\\
F_{n,\eps}(\xi)=F_{\eps}(\xi)\ \ \mbox{if}\ \ \xi\in
B_{\rho_n+1}(0),\qquad\qquad\label{part3.3.5.2}\\
F_{n,\eps}(\xi)=\nu_n(1+|\xi|^2)\ \ \mbox{if}\ \ |\xi|\geq
\rho_n+2,\qquad\label{part3.3.5.3}\\
|F_{n,\eps}(\xi)|\leq \nu_n(1+|\xi|^2)\ \ \mbox{for all}\ \
\xi\in\rr^N,\qquad\label{part3.3.5.4}
\end{eqnarray}
where $\nu_n$ is a positive constant which depends only on $\rho_n$ and $p$.\\
With $F_{n,\eps}$ defined above we consider the problem:
\begin{equation}\label{part3.3.6}
\left\{\begin{array}{l}
\di\frac{\partial u}{\partial t}-\Delta u=F_{n,\eps}(\nabla
u)\ \mbox{in}\ Q_T,\\
\di\frac{\partial u}{\partial\nu}=0\ \ \ \ \mbox{on}\ \ \Gamma_T,\\
u(0,\cdot)=u_0^n\ ~\mbox{in}\ ~\Omega.\end{array}\right.
\end{equation}
Thanks to the regularity of $u_0^n$ and to relations \eq{part3.3.1}, \eq{part3.3.5.5},
\eq{part3.3.5.2}, \eq{part3.3.5.3} and \eq{part3.3.5.4} we can apply Theorem
V.7.4 in \cite{LSU} to the problem \eq{part3.3.6}. Thus, there exists $u^{n,\eps}\in
C^{1+\alpha/2,2+\alpha}(\overline{Q_T})$,$\alpha\in(0,1)$, the unique
solution of problem \eq{part3.3.6}.
For any $(t,x)\in Q_T$ let denote by:
\begin{equation}\label{part3.3.6.1}
f_{n,\eps}(t,x)=F_{n,\eps}(\nabla u^{n,\eps})(t,x).
\end{equation}
Then, thanks to the regularity of $F_{n,\eps}$ and $u^{n,\eps}$ it follows that:
$f_{n,\eps}\in C^{\frac{1+\alpha}{2},{1+\alpha}}(\overline{Q_T})$
and  $u^{n,\eps} $ verifies:
\begin{equation}\label{part3.3.7}
\left\{\begin{array}{l}
\di\frac{\partial u^{n,\eps}}{\partial t}-\Delta u^{n,\eps}=f_{n,\eps}
  \ \ \mbox{in}\ \ Q_T,\\
\di\frac{\partial u^{n,\eps}}{\partial\nu}=0\ \ ~\mbox{on}\ \ ~\Gamma_T,\\
u^{n,\eps}(0,\cdot)=u^{n}_0\ \ \mbox{in}\ \ ~\Omega.
\end{array}\right.
\end{equation}
Applying Theorem III.12.2 in
\cite{LSU}, on the local regularity of solution for parabolic problem of
\eq{part3.3.7} type, we get:
$$u^{n,\eps}\in C_{loc}^{\frac{3+\alpha}{2},
  3+\alpha}(Q_T)\cap C^{1+\alpha/2, 2+\alpha}(\overline{Q_T}).$$
In the sequel, we show that, for $\eps\in(0,1)$,
\begin{equation}\label{part3.3.7.3}
|\nabla u^{n,\eps}(t,x)|\leq \rho_n,\ \ \forall\ (t,x)\in Q_T,
\end{equation}
where $\rho_n$ is given by \eq{part3.3.5.1}.
For this, we will use the Bernstein technique. First we introduce the
parabolic operator ${\mathcal L}$  defined on $C^{0,1}(\overline{Q_T})\cap
C^{1,2}(Q_T)$ by:
$${\mathcal L}(v)=\frac{\partial v}{\partial t}-\Delta v+b(t,x)\cdot\nabla v,$$
where $b\in [L^{\infty}(Q_T)]^N$ is given by:
$$b(t,x)=-(\nabla F_{n,\eps})(\nabla u^{n,\eps})(t,x)\ \mbox{in}\ \ Q_T.$$
Setting $w=|\nabla u^{n,\eps}|^2$, then $w\in C^{0,1}(\overline{Q_T})\cap
C^{1,2}(Q_T)$  and verifies:
$${\mathcal
  L}(w)=-2\sum\limits_{i,j=1}^N\left(\frac{\partial^2u^{n,\eps}}{\partial
    x_i\partial x_j}\right)^2\leq 0,$$
hence, thanks to Lemma \pr{part3.l.3} and to relations \eq{part3.3.1} and
\eq{part3.3.5.1} we have:
$$
\di\frac{\partial w}{\partial\nu}\leq 0\ \ \mbox{on}\ \ \Gamma_T\ \
\mbox{and}\ \ \
\di w(0,x)\leq \rho_n^2\ \ \mbox{in}\ \ \ \Omega.$$
Then, by the Comparison Principle (Lemma \pr{part3.l.4}), we obtain:
$w\leq \rho_n^2\ \ \mbox{in}\ \ Q_T$
and relation \eq{part3.3.7.3} is proved.\\
Combining \eq{part3.3.5.2}, \eq{part3.3.6} and \eq{part3.3.7.3}, we finally
obtain that:
$u^{n,\eps}\in C_{loc}^{\frac{3+\alpha}{2},
  3+\alpha}(Q_T)\cap C^{1+\alpha/2, 2+\alpha}(\overline{Q_T})$ is the
solution of the initial boundary value problem:
 \begin{equation}\label{part3.3.7.4}
\left\{\begin{array}{l}
\di\frac{\partial u}{\partial t}-\Delta u=F_{\eps}(\nabla
u)\ \mbox{in}\ Q_T,\\
\di\frac{\partial u}{\partial\nu}=0\ \ \ \ \mbox{on}\ \ \Gamma_T,\\
u(0,\cdot)=u_0^n\ ~\mbox{in}\ ~\Omega.\end{array}\right.
\end{equation}
Moreover, we notice that, in  \eq{part3.3.7.4}, $F_{\eps}$ is independent
of $n$. \\

{\bf Second step: ``Estimates for $\ u^{n,\eps}\ $''.}\\
For $\eps>0$ and  $n$ a positive entire let set:
\begin{equation}\label{part3.3.7.4.1}
m_{n,\eps}=m(0)+\frac{1}{n}-a\eps^{p/2}\cdot T
\end{equation}
and
\begin{equation}\label{part3.3.7.4.2}
M_{n,\eps}=M(0)+\frac{2}{n}+a\eps^{p/2}\cdot T.
\end{equation}
The next proposition gives some estimates of $u^{n,\eps}$ which will
allow us to pass to the limits in  \eq{part3.3.7.4}, as $\eps$ tends to $0$:
\begin{prop}\label{part3.p.1}
For all $p\in(0,+\infty)$, the solution \\
$u^{n,\eps}\in C_{loc}^{\frac{3+\alpha}{2},
  3+\alpha}(Q_T)\cap C^{1+\alpha/2, 2+\alpha}(\overline{Q_T})$ of problem \eq{part3.3.7.4}  satisfies:
\begin{equation}\label{part3.3.8}
m_{n,\eps}\leq u^{n,\eps}\leq M_{n,\eps}\ \ \mbox{in}\ \ Q_T,
\end{equation}
\begin{equation}\label{part3.3.8.1}
\|\nabla
u^{n,\eps}(t)\|_{\infty}\leq\left(\demi\right)^{1/2}(M_{n,\eps}-m_{n,\eps}+\frac{1}{n})\cdot
t^{-\demi}\ \mbox{for all}\  t\in(0,T),
\end{equation}
and, if $p\neq 1$:
\begin{equation}\label{part3.3.9}
\|\nabla
u^{n,\eps}(t)\|_{\infty}\leq\left(\frac{\max\{p,2\}}{ap|1-p|}\right)^{1/p}(M_{n,\eps}-m_{n,\eps}+\frac{1}{n})^{1/p}\cdot
  t^{-1/p}\ \mbox{for all}\  t\in(0,T).
\end{equation}
\end{prop}
\proof
The two inequalities in \eq{part3.3.8} are simple consequences of Lemma
\pr{part3.l.4}. \\
Instead, to prove  \eq{part3.3.8.1} and \eq{part3.3.9} we will use the
Bernstein technique and the proof is similar to that given in \cite{BL},\cite{ GGK} and
\cite{Li2}.
Let denote by $w$ the function defined on $Q_T$ by:
\begin{equation}\label{part3.3.10}
w=\frac{|\nabla u^{n,\eps}|^2}{\theta(u^{n,\eps})}.
\end{equation}
where $\theta$ is a strict positive function of
$C^2([m_{n,\eps},M_{n,\eps}])$ class, which will be chosen later according
to the exponent $p$.\\
Then, thanks to the regularity of function $u^{n,\eps}$ we have:
$$w\in C^{0,1}(\overline{Q_T})\cap C^{1,2}(Q_T).$$
Moreover:
$$
\begin{array}{l}
\di\frac{\partial
w}{\partial\nu}=\di\frac{1}{[\theta(u^{n,\eps})]^2}\left(\di\frac{\partial
|\nabla
  u^{n,\eps}|^2}{\partial\nu}\theta(u^{n,\eps})-|\nabla
u^{n,\eps}|^2\cdot\theta'(u^{n,\eps})\di\frac{\partial
  u^{n,\eps}}{\partial\nu}\right)\\
  \qquad=\di\frac{1}{[\theta(u^{n,\eps})]}\cdot\frac{\partial
|\nabla u^{n,\eps}|^2}{\partial\nu}\ \ \mbox{on}\ \Gamma_T.
\end{array}$$ Since $\theta$ is a positive function, this last
relation and Lemma \pr{part3.l.3} imply:
\begin{equation}\label{part3.3.10.3}
\frac{\partial w}{\partial\nu}\leq 0\ \mbox{on}\ \Gamma_T.
\end{equation}
Denote by
${\mathcal N}$ the semi-linear parabolic operator defined on
$C^{0,1}(\overline{Q_T})\cap C^{1,2}(Q_T)$ by:
$${\mathcal N}(v)=\frac{\partial v}{\partial t}-\Delta v-b(t,x)\cdot\nabla
v-c(t,x)v^2-d(t,x)v^{1+p/2}$$
where
$$b(t,x)=(\nabla F_{\eps})(\nabla u^{n,\eps})(t,x)+\frac{2\theta'(u^{n,\eps})\nabla u^{n,\eps}(t,x)}{\theta(u^{n,\eps})(t,x)},$$
$$c(t,x)=\theta''(u^{n,\eps})(t,x),$$
and
\begin{equation}\label{part3.3.10.9}
d(t,x)=a(p-1)\theta^{\frac{p-2}{2}}(u^{n,\eps})(t,x)
\theta'(u^{n,\eps})(t,x).
\end{equation}
The function $w$ being introduced by \eq{part3.3.10} we have:
\begin{equation}\label{part3.3.10.1}
\begin{array}{l}\di {\mathcal N}(w)=-\frac{2}{\theta(u^{n,\eps})}\di\sum_{i,j=1}^{N}\left(\frac{\partial
  ^{2} u^{n,\eps}}{\partial x_i\partial x_j}\right)^2+\frac{\theta'(u^{n,\eps})}{\theta^2(u^{n,\eps})}[(\nabla F_{\eps})(\nabla u^{n,\eps})\cdot\nabla
u^{n,\eps}\\
\qquad\qquad\qquad\qquad \ \ \ \ \ \ \ \ \ \ \ \ \ \ \ \ \ \ -F_{\eps}(\nabla u^{n,\eps})-a(p-1)|\nabla u^{n,\eps}|^p]|\nabla u^{n,\eps}|^2.\end{array}
\end{equation}
To prove \eq{part3.3.8.1} we will distinguish between the two cases below.\\
$i)$ The case $0<p\leq 1$. We take $\theta$  in \eq{part3.3.10} as follows:
$$\theta(\xi)=\demi(M_{n,\eps}-m_{n,\eps}+\frac{1}{n})^2-\demi(M_{n,\eps}-\xi)^2,\
\ \
\xi\in[m_{n,\eps},M_{n,\eps}],$$
where $m_{n,\eps}$ and $M_{n,\eps}$ are defined by \eq{part3.3.7.4.1} and \eq{part3.3.7.4.2}.
So $\theta$ verifies:
$$\theta(\xi)\geq\frac{1}{2n^2},\ \ \ \ \ \ \ \theta
'(\xi)=M_{n,\eps}-\xi,\ \ \ \ \ \ \ \ \ \
\theta''(u)=-1,$$
and we deduce that:
$$\theta'(u^{n,\eps})\geq 0$$
and
$$d(t,x)=a(p-1)\theta^{\frac{p-2}{2}}(u^{n,\eps})(t,x) \theta'(u^{n,\eps})(t,x)\leq 0.$$
Combining these last points with \eq{part3.3.4} and \eq{part3.3.10.1}  it
follows that:
\begin{equation}\label{part3.3.10.2}
{\mathcal N}(w)\leq 0.
\end{equation}
Taking into account \eq{part3.3.5.1} and \eq{part3.3.7} we have:
$$w(0)=\frac{|\nabla u_0^n|^2}{\theta(u_0^n)}\leq 2\rho_n^2n^2.$$
So, for  $n$ a fixed entire, choose $\eta>0$ such that:
\begin{equation}\label{part3.3.10.4}
w(0)\leq 2\rho_n^2n^2\leq\frac{1}{\eta},
\end{equation}
and denote by $v$ the function defined on $Q_T$ by:
$$v(t,x)=(t+\eta)^{-1}.$$
Since $a>0$ and $p\in(0,1)$ we have:
\begin{equation}\label{part3.3.10.5}
{\mathcal N}(v)=-d(t,x)\cdot(t+\eta)^{-(1+p/2)}\geq 0.
\end{equation}
So, recalling \eq{part3.3.10.3},\eq{part3.3.10.2},
\eq{part3.3.10.4}, \eq{part3.3.10.5} and  Lemma \pr{part3.l.4} we
get:
$$w(t,x)\leq (t+\eta)^{-1}<t^{-1} \ \ \mbox{for all}\ (t,x)\in Q_T,$$
and we deduce that \eq{part3.3.8.1} holds for $p\in(0,1]$.\\

$ii)$ The case $p>1$. In \eq{part3.3.10} we consider the function $\theta$
defined by:
$$\theta(\xi)=\demi(M_{n,\eps}-m_{n,\eps}+\frac{1}{n})^2-\demi(\xi-m_{n,\eps})^2,\
\ \ \xi\in[m_{n,\eps},M_{n,\eps}],$$ then $\theta$ satisfies:
$$\theta(\xi)\geq\frac{1}{2n^2},\ \ \ \ \ \ \ \theta '(\xi)=m_{n,\eps}-\xi,\ \
\ \ \ \ \ \ \ \ \ \theta''(\xi)=-1,$$ and we deduce that:
$$\theta'(u^{n,\eps})\leq 0,$$
and
 $$d(t,x)=a(p-1)\theta^{\frac{p-2}{2}}(u^{n,\eps})(t,x) \theta'(u^{n,\eps})(t,x)\leq 0.$$
Combining these last points with \eq{part3.3.5} and \eq{part3.3.10.1}  it
follows that:
$${\mathcal N}(w)\leq 0.$$
As previously, we can prove \eq{part3.3.8.1} for the case $p\geq 1$ by comparing $w$ and $v$.\\

To prove \eq{part3.3.9} we will distinguish among three cases:\\
$i)$ The case $0<p< 1$. In \eq{part3.3.10}, we consider the following function:
$$\theta(\xi)=(\frac{2}{ap(1-p)})^{2/p}(M_{n,\eps}-m_{n,\eps}+\frac{1}{n})^{\frac{2-p}{p}}\cdot(\xi-m_{n,\eps}+\frac{1}{n}),
\ \xi\in[m_{n,\eps}, M_{n,\eps}].$$
Thus
\begin{equation}\label{part3.3.10.6}
\theta(\xi)\geq \left[\frac{2}{anp(1-p)}\right]^{2/p},\ \ \
\theta'(\xi)\geq 0\ \ \ \mbox{and}\ \ \ \theta''(\xi)=0,\ \xi\in[m_{n,\eps}, M_{n,\eps}].
\end{equation}
$w$ being given by \eq{part3.3.10},  thanks to relations \eq{part3.3.4} and
\eq{part3.3.10.1} we obtain:
\begin{equation}\label{part3.3.12}
{\mathcal N}(w)\leq 0
\end{equation}
Taking into account \eq{part3.3.6} and  \eq{part3.3.10.6}, we can choose
$\eta>0$ such that:
\begin{equation}\label{part3.3.14}
w(0)=\frac{|\nabla u_0^n|^2}{\theta(u_0^n)}\leq 2\rho_n^2\left[\frac{anp(1-p)}{2}\right]^{2/p}<\frac{1}{\eta^{2/p}},
\end{equation}
where $\rho_n$ is given by \eq{part3.3.5.1}.\\
Let $v$ be a function defined on $Q_T$ by:
$$v(t,x)=(t+\eta)^{-2/p}.$$
With $d$ given by \eq{part3.3.10.9}, and $\theta$ being chosen as above we have:
$$d(t,x)=-\frac{2}{p}\left(\frac{M_{n,\eps}-m_{n,\eps}+\frac{1}{n}}{u-m_{n,\eps}+\frac{1}{n}}\right)^{\frac{2-p}{2}}\leq
-\frac{2}{p}$$
And we deduce that:
\begin{equation}\label{part3.3.15}
{\mathcal N}(v)=(-d-2/p)(t+\eta)^{-\frac{p+2}{2}}\geq 0
\end{equation}
Combining relations
\eq{part3.3.10.3},\eq{part3.3.12},\eq{part3.3.14},\eq{part3.3.15},
and Lemma \pr{part3.l.4} we get:
$$w(t,x)\leq (t+\eta)^{2/p}<t^{-2/p}\ \mbox{for all}\ t>0,$$
and we deduce \eq{part3.3.9}, for $0<p\leq 1$.\\
$ii)$ The case $1<p<2$. In \eq{part3.3.10}, we choose the
following function $\theta$:
$$\theta(\xi)=\Big[\frac{2}{ap(p-1)}\Big]^{2/p}(M_{n,\eps}-m_{n,\eps}+\frac{1}{n})^{\frac{2-p}{p}}(M_{n,\eps}-\xi+\frac{1}{n}),\
\xi\in[m_{n,\eps}, M_{n,\eps}]$$
Thus:
\begin{equation}\label{part3.3.10.10}
\theta(\xi)\geq \Big[\frac{2}{anp(p-1)}\Big]^{2/p},\ \ \theta'(\xi)\leq
0\ \ \mbox{and}\ \  \theta''(\xi)=0,\ \xi\in[m_{n,\eps}, M_{n,\eps}].
\end{equation}
and we get \eq{part3.3.9} as previously.\\
$iii)$ The case $p\geq 2$. This time we prove \eq{part3.3.9} in the two
cases above, by taking:
$$\theta(\xi)=\Big[\frac{1}{a(p-1)}(M_{n,\eps}-\xi+\frac{1}{n})\Big]^{2/p},\ \xi\in[m_{n,\eps}, M_{n,\eps}].$$
\endproof
We came back to the problem \eq{part3.3.7.4} and we notice that $$\eps\mapsto
F_{\eps}(\xi)\ \mbox{is a nondecreasing function for}\ \ 0<p\leq 1$$
and
$$\eps\to F_{\eps}(\xi)\ \mbox{is a decreasing function for}\ \ p> 1.$$
Then, thanks to relations \eq{part3.3.3} and  \eq{part3.3.8} we
can apply Lemma \pr{part3.l.4} and we obtain that the set
$(u^{n,\eps})_{\eps>0}$ is bounded and monotone with respect to
$\eps$, and consequently, there exists $u^{n}\in L^{\infty}(Q_T)$
such that
$$u^{n,\eps}\nearrow u^{n}\ \mbox{in}\ Q_T\ ~\mbox{as}\
\eps\searrow 0,\ \ \mbox{if}\ \ 0<p\leq 1$$
and
$$u^{n,\eps}\searrow u^{n}\ \mbox{in}\ Q_T\ ~\mbox{as}\
\eps\searrow 0,\ \ \mbox{if}\ \ p> 1.$$ Moreover, from relations
\eq{part3.3.1} and \eq{part3.3.8}, the hypotheses of Theorem V.7.2
in \cite{LSU} are satisfied and we deduce that the solutions
$u^{n,\eps}$ of \eq{part3.3.7.4} verify:
\begin{equation}\label{part3.3.16}
\|u^{n,\eps}\|_{C^{\frac{1+\delta}{2},1+\delta}(Q_T)}\leq C
\end{equation}
where $\delta\in(0,1)$ and $C$ are two positive constants which depend only
on $m,M,\|u_0^n\|_{\Omega}^{(2)}$ and $\Omega$. Thus, we deduce that for
all $n,$ the set
$\{u^{n,\eps}, 0<\eps<1\}$ is bounded in
$C^{\frac{1+\delta}{2},1+\delta}(\overline{Q_T})$.
Let be $f_{n,\eps}$ the function given by \eq{part3.3.6.1}, then, thanks to
 the regularity of $F_{\eps}$ and to \eq{part3.3.16},
 the set $\{f_{n,\eps}, 0<\eps<1\}$ is bounded in
$C^{\delta/2,\delta}(\overline{Q_T})$. Since $u^{n,\eps}\in
C^{\frac{1+\delta}{2},1+\delta}(\overline{Q_T}) $ is the solution
of problem \eq{part3.3.7.4}, the hypotheses of Theorem IV.5.3 in
\cite{LSU} on the regularity in H\"older spaces of solutions for
parabolic equations, are verified and therefore we get:
$$u^{n,\eps}\in C^{1+\delta/2,2+\delta}(\overline{Q_T}),$$
moreover, there exists a constant $C>0$,  not depending on
$\eps\in(0,1)$, such that:
\begin{equation}\label{part3.3.17}
\|u^{n,\eps}\|_{C^{1+\delta/2,2+\delta}(\overline{Q_T})}\leq C(\|u_0^n\|_{C^{2+\delta}(\overline{\Omega})}+\|f_{n,\eps}\|_{C^{\delta/2,\delta}(\overline{Q_T})}\leq
C_n
\end{equation}
Thus, the set $\{u^{n,\eps}, 0<\eps<1\}$ is bounded in
$C^{1+\delta/2,2+\delta}(\overline{Q_T})$. Since for any $0\leq\nu<\delta$
$$C^{1+\delta/2,2+\delta}(\overline{Q_T})\hookrightarrow C^{1+\nu/2,2+\nu}(\overline{Q_T}),$$
with compact embedding,
we deduce that $\{u^{n,\eps}, 0<\eps<1\}$ is a precompact set in
$C^{1+\nu/2,2+\nu}(\overline{Q_T})$ and it follows that, ``to a
subsequence'' we have:
\begin{equation}\label{part3.3.18}
u^{n,\eps}\to
u^{n}\ ~\mbox{in}\ C^{1+\nu/2,2+\nu}(\overline{Q_T})\ \mbox{as}\ \eps\searrow
0
\end{equation}
On the other hand, for all $\xi\in\rr^N$:
$$F_{\eps}(\xi)\to a|\xi|^p\ \mbox{as}\ \eps\searrow 0,$$
So, we can pass to the limit in  \eq{part3.3.7.4}, as
$\eps\searrow 0$, and we obtain that
$u^{n}\in C^{1+\nu/2,2+\nu}(\overline{Q_T})$ is a solution of the following
initial boundary value problem:
\begin{equation}\label{part3.3.19}
\left\{\begin{array}{l}
\di\frac{\partial u^{n}}{\partial t}-\Delta u^n=a|\nabla
u^n|^p\ \ ~\mbox{in}\ \ Q_T,\\
\di\frac{\partial u^{n}}{\partial\nu}=0\ \ ~\mbox{on}\ \ \Gamma_T,\\
u^{n}(0,x)=u_0^n(x)\ \ ~\mbox{in}\ \ \Omega.\end{array}\right.
\end{equation}
Applying the Comparison Principle, [Theorem 1 in \cite{Kp}],
 we get also that this solution is unique in $C^{1,2}(\overline{Q_T})$. \\

{\bf Third step: ``Estimates for $u^n\ $''.}\\
The aim of the following proposition is to prove that $(u^n)_n$
satisfies also the estimates \eq{part3.3.8}, \eq{part3.3.8.1} and \eq{part3.3.9} for
$\eps=0$, and is bounded in a H\"older space.
\begin{prop}\label{part3.p.2}
The solution $u^n\in C^{1+\nu/2,2+\nu}(\overline{Q_T})$ of problem
\eq{part3.3.19} satisfies the following properties:
\begin{equation}\label{part3.3.20}
m(0)+\frac{1}{n}\leq u^n(t,x)\leq M(0)+\frac{2}{n},
\end{equation}
\begin{equation}\label{part3.3.20.1}
\|\nabla u^n(t)\|_{\infty}\leq\left(\demi\right)^{1/2}(M(0)-m(0)+\frac{2}{n})\cdot
t^{-\demi},\ \mbox{for all}\  t\in(0,T),
\end{equation}
and, if $p\neq 1$ then:
\begin{equation}\label{part3.3.21}
\|\nabla
u^n(t)\|_{\infty}\leq\left(\frac{\max\{p,2\}}{ap|1-p|}\right)^{1/p}(M(0)-m(0)+\frac{2}{n})^{1/p}\cdot
t^{-1/p},\ \mbox{for all}\  t\in(0,T).
\end{equation}
Moreover, there exists $\delta\in(0,1)$ such that, for all $\tau\in(0,T)$:
\begin{equation}\label{part3.3.22}
\mbox{the sequence} \ (u^n)_n\ \mbox{is bounded in}\ \  C^{1+\delta/2,2+\delta}(\overline{Q_{\tau,T}}).
\end{equation}
(This bound depends only on
   $\tau,\Omega,p,m(0)$ and $M(0)$.)
\end{prop}
\proof
Relations \eq{part3.3.20}, \eq{part3.3.20.1} and \eq{part3.3.21} are direct
consequences of \eq{part3.3.8}, \eq{part3.3.8.1}, \eq{part3.3.9} and
\eq{part3.3.18}. In order to prove \eq{part3.3.22} we denote by $f_n$ the
function defined on $Q_T$ by:
$$f_n(t,x)=a|\nabla u^n|^p(t,x).$$
Then $u^n\in C^{1+\frac{\nu}{2},2+\nu}(\overline{Q_T})$ is the solution of
the following problem:
\begin{equation}\label{part3.3.24}
\left\{\begin{array}{l}
\di\frac{\partial u^{n}}{\partial t}-\Delta u^n=f_n\ ~\mbox{in}\ Q_T,\\
\di\frac{\partial u^{n}}{\partial\nu}=0\ \ ~\mbox{on}\ \Gamma_T,\\
u^{n}(0,\cdot)=u_0^n\ ~\mbox{in}\ \Omega.\end{array}\right.
\end{equation}
Consider $\tau\in(0,T)$. Thanks to relation \eq{part3.3.20.1}, $f_n\in L^{\infty}(Q_{\tau,T})$
and:
\begin{equation}\label{part3.3.23}
\|f_n\|_{L^{\infty}(Q_{\tau,T})}\leq\frac{a}{2^{p/2}}(M(0)-m(0)+2)^p\tau^{-p/2},\
  \ \ \forall \ n\in\nn
\end{equation}
Consequently, the sequence $(f_n)_{n\geq 0}$ is uniformly bounded in $L^{\infty}(Q_{\tau,T})$.\\
In the sequel, we decompose the problem \eq{part3.3.24}  into two parts. \\
On the one hand, we denote by $v^n$ the solution of the heat equation on $Q_{\tau/3,T}$:
\begin{equation}\label{part3.3.25}
\left\{\begin{array}{l}
\di\frac{\partial v^{n}}{\partial t}-\Delta v^n=0\ \ ~\mbox{in}\ \ Q_{\tau/3,T},\\
\di\frac{\partial v^{n}}{\partial\nu}=0\ \ ~\mbox{on}\ \ \Gamma_{\tau/3,T},\\
v^{n}(\tau/3,x)=u^n(\tau/3,x)\ \ ~\mbox{in}\ \ \Omega.\end{array}\right.
\end{equation}
Thanks to the regularity effect of the heat equation it follows that:
\begin{equation}\label{part3.3.61}
v^n\in C^{\infty}(\overline{Q_{2\tau/3,T}})
\end{equation}
and from Lemma \pr{part3.l.5} and relations \eq{part3.2.8} and \eq{part3.3.20},
for all $l\in\nn$ and $\alpha\in\nn^N$ we have:
\begin{equation}\label{part3.3.26.1}
\|D_x^{\alpha}D_t^l v^n\|_{\infty,Q_{2\tau/3,T}}\leq
C(T,\Omega)(M(0)+m(0)+1)\tau^{-(\frac{|\alpha|}{2}+l)}.
\end{equation}
Next, we denote by $w^n$ the solution of the problem:
\begin{equation}\label{part3.3.26}
\left\{\begin{array}{l}
\di\frac{\partial w^{n}}{\partial t}-\Delta w^n=f_n(t,x)\ ~\mbox{in}\ Q_{\tau/3,T},\\
\di\frac{\partial w^{n}}{\partial\nu}=0\ \ ~\mbox{on}\ \Gamma_{\tau/3,T},\\
w^{n}(\tau/3,\cdot)=0\ ~\mbox{in}\ \Omega.\end{array}\right.
\end{equation}
Taking into account \eq{part3.3.18}, we have $f_n\in
C^{\frac{1+\nu}{2},1+\nu}(\overline{Q_T})$ and we deduce that $w^n\in C^{1+\nu/2,2+\nu}(\overline{Q_{\tau/3,T}})$.
Since $f_n\in L^{\infty}(Q_{\tau/3,T})$, we have in particular $f_n\in
L^q(Q_{\tau/3,T})$ for all
$q>1$. Thus, we can apply Theorem 7.20 in \cite{L}, on the regularity of
parabolic solutions in
$L^q$ spaces, and we get that, there exists a constant $C>0$, independent on $n$,
such that:
\begin{equation}\label{part3.3.27}
\|D_x^2 w^n\|_{q,Q_{\frac{\tau}{3},T}}+\|D_t w^n\|_{q,Q_{\frac{\tau}{3},T}}\leq
C\|f_n\|_{q,Q_{\frac{\tau}{3},T}}\leq C|Q_{\frac{\tau}{3},T}|^{1/q}\|f_n\|_{\infty,Q_{\frac{\tau}{3},T}}.
\end{equation}
Combining \eq{part3.3.23} with \eq{part3.3.27} we get:
\begin{equation}\label{part3.3.28}
\|D_t w^n\|_{q,Q_{\frac{\tau}{3},T}}+\|D_x^2 w^n\|_{q,Q_{\frac{\tau}{3},T}}\leq C(M(0),m(0),p,\tau, T, \Omega).
\end{equation}
Since $u^n=v^n+w^n$, from \eq{part3.3.26.1} and \eq{part3.3.28} we get on
the one hand:
\begin{equation}\label{part3.3.29}
\|D_t u^n\|_{q,Q_{2\tau/3,T}}+\|D_x^2 u^n\|_{q,Q_{2\tau/3,T}}\leq C(M(0),m(0),p,q,\tau, T, \Omega).
\end{equation}
On the other hand relations \eq{part3.3.20} and \eq{part3.3.20.1} yield:
\begin{equation}\label{part3.3.30}
\|u^n\|_{\infty,Q_{2\tau/3,T}}\leq C_1(M(0),m(0))
\end{equation}
and
\begin{equation}\label{part3.3.31}
\|D_x u^n\|_{\infty,Q_{2\tau/3,T}}\leq C_2(M(0),m(0),p,\tau).
\end{equation}
So, combining \eq{part3.3.29}, \eq{part3.3.30} and \eq{part3.3.31} we get:
$$\|u^n\|_{W_q^{1,2}(Q_{2\tau/3,T})}\leq C(M(0),m(0),p,q,\tau,T, \Omega)\
\mbox{for all}\ n\in\nn\ \mbox{and}\ q>1.$$
We choose $q>N+2$, then, applying Lemma II.3.3 in \cite{LSU} (on the
embedding of Sobolev spaces into H\"older spaces),  we deduce that, for any
$\beta$ satisfying $0<\beta<1-\frac{N+2}{q}$, there exists a constant $C>0$
such that:
$$\|\nabla u^n\|_{C^{\beta/2,\beta}(\overline{Q_{2\tau/3,T}})}\leq
C(q,\beta,\tau,T,\Omega)\|u^n\|_{W_q^{1,2}(Q_{2\tau/3,T})}.$$
Since the sequence $(u_n)_n$ is bounded in
$W_q^{1,2}(Q_{2\tau/3,T})$, we deduce that $(|\nabla u_n|)_n$ is bounded in $C^{\beta/2,\beta}(\overline{Q_{2\tau/3,T}})$.
Consequently the sequence $(f_n=a|\nabla u_n|^p)_n$ is uniformly bounded in
$C^{\delta/2,\delta}(Q_{2\tau/3,T})$, where  $\delta=\delta(\beta,p)$.

We came back  to problems \eq{part3.3.24}, \eq{part3.3.25} and
\eq{part3.3.26} in $Q_{2\tau/3,T}$. By reiterating the process above we
get, thanks to Theorem IV.5.3 in \cite{LSU}, that:\\
$i)$ $w^n\in C^{1+\delta/2,2+\delta}(\overline{Q_{2\tau/3,T}})$ and there
exists a constant $C>0$, independent on $n$
such that:
\begin{equation}\label{part3.3.32}
\|w^n\|_{C^{1+\delta/2,2+\delta}(\overline{Q_{\frac{2\tau}{3},T}})}\leq
  C\|f_n||_{C^{\delta/2,\delta}(\overline{Q_{\frac{2\tau}{3},T}})}\leq C(m(0),M(0),p,N, \tau,T,\Omega).
\end{equation}
$ii)$ $v^n$ satisfies relation \eq{part3.3.26.1} on $Q_{\tau,T}$.\\

Thus, recalling  \eq{part3.3.26.1} and \eq{part3.3.32} we obtain that $u^n\in
C^{1+\delta/2,2+\delta}(\overline{Q_{\tau,T}})$ and:
\begin{equation}\label{part3.3.33}
\|u^n\|_{C^{1+\delta/2,2+\delta}(\overline{Q_{2\tau/3,T}})}\leq C(m(0),M(0),p,N,\tau,T,\Omega),
\end{equation}
which ends the proof of Proposition \pr{part3.p.2}.
\endproof

{\bf Four step: ``Proof of the existence of solutions''.}\\
On the one hand, thanks to the Comparison Principle, [Theorem 1 in
\cite{Kp}], and to relations \eq{part3.3.1} and \eq{part3.3.20} the
sequence $(u^n)_n$ is decreasing and uniformly bounded. Consequently, there
exists  $u\in L^{\infty}(Q_T)$ such that:
\begin{equation}\label{part3.3.34}
u^n\searrow u\ \ ~\mbox{in}\ \ ~Q_T.\ \
\end{equation}
On the other hand, by Proposition \pr{part3.p.2} we deduce that, for any
$\tau\in(0,T)$,  the sequence $(u^n)_{n\geq 1}$ is bounded in
$C^{1+\delta/2,2+\delta}(\overline{Q_{\tau,T}})$. Since for all $\nu\in(0,\delta)$:
$$\ C^{1+\delta/2,2+\delta}(\overline{Q_{\tau,T}})\hookrightarrow
C^{1+\nu/2,2+\nu}(\overline{Q_{\tau,T}})$$
with compact embedding,
``to a subsequence'', we have:
\begin{equation}\label{part3.3.35}
u^n\to u\ \ \mbox{in}\ \ C^{1+\nu/2,2+\nu}(\overline{Q_{\tau,T}})\ \
\mbox{as}\ \ n\to\infty.
\end{equation}
Hence, $u\in C^{1+\nu/2,2+\nu}(\overline{Q_{\tau,T}})$ and thanks
to relations \eq{part3.3.34} and \eq{part3.3.35} we may let
$t\to\infty$ in the first and the second equation of problem
\eq{part3.3.19} and we obtain that, for all $\tau\in(0,T)$, $u$
satisfies:
\begin{equation}\label{part3.3.35.1}
\left\{\begin{array}{l}
\di\frac{\partial u}{\partial t}-\Delta u=a|\nabla
u|^p\ ~\mbox{in}\ Q_{\tau,T},\\
\di\frac{\partial u}{\partial\nu}=0\ \ ~\mbox{on}\ \Gamma_{\tau,T}.\end{array}\right.
\end{equation}
 Moreover, passing to limits in \eq{part3.3.20.1} and  \eq{part3.3.21}, as
 $n$ tends to $\infty$, we get \eq{part3.1.5.1} and \eq{part3.1.6}. The
 relations \eq{part3.1.4} and
\eq{part3.1.5} are direct consequences of Lemma \pr{part3.l.4}.\\
 So, we have to identify the initial data $\mu_0$.
 For
 $t\in(0,T)$, let denote by
 $v(t)=S(t)\mu_0$ and $v^n=S(t)u_0^n$, where  $(S(t))_{t\geq 0}$ is the
 heat semigroup in $\lq, q\geq 1$, for the homogeneous Neumann boundary
 value problem. Then, by the Comparison Principle [Lemma \pr{part3.l.4}],
 for  $n\in\nn$ we have:
$$v^n\leq u^n\ \ \mbox{in}\ \ Q_T.$$
Using \eq{part3.3.34}, we may let  $n\to\infty$ in the above inequality and
we obtain:
\begin{equation}\label{part3.3.36}
v\leq u\ \ \mbox{in}\ \ Q_T.
\end{equation}
Since $v\in C(\overline{Q_T})$,   it follows that:
\begin{equation}\label{part3.3.37}
\mu_0(x_0)=\lim\limits_{\stackrel{(t,x)\to(0,x_0)}{(t,x)\in Q_T}}v(t,x)\leq \liminf\limits_{\stackrel{(t,x)\to(0,x_0)}{(t,x)\in Q_T}}u(t,x),
\end{equation}
for any  $x_0\in\Omega$. Furthermore for  $n\in\nn$ we have:
$$u\leq u^n\ \ \mbox{in}\ \ Q_T,$$
then:
$$\limsup\limits_{\stackrel{(t,x)\to(0,x_0)}{ (t,x)\in
  Q_T}}u(t,x)\leq\limsup\limits_{\stackrel{(t,x)\to(0,x_0)}{ (t,x)\in Q_T}}u^n(t,x)=
u_0^n(x_0).$$
Since $(u_0^n)_n$ is a decreasing sequence and converges to $\mu_0$ we can
pass to the limits in the above inequality and we get
\begin{equation}\label{part3.3.38}
\limsup\limits_{\stackrel{(t,x)\to(0,x_0)}{ (t,x)\in
  Q_T}}u(t,x)\leq\mu_0(x_0).
\end{equation}
Combining \eq{part3.3.37}, \eq{part3.3.38} and the fact that  $x_0$ is
anywhere in $\Omega$ we deduce that $u\in C(\overline{Q_T})\cap
C^{1+\nu/2,2+\nu}(\overline{Q_{\tau,T}})$ is a classical solution of the
problem \eq{part3.1.1}. Which end the existence proof of solutions of
problem \eq{part3.1.1}, for $a>0$.\\

{\bf Fifth step: ``Uniqueness of the solution''.}\\
The uniqueness is a direct consequence of the following lemma:
\begin{lemma}\label{part3.l.6}
Let $a>0$, $p>1$, $\Omega\subset\rr^N$ a bounded and convex open set with
smooth boundary. Let  $\mu_0\in C(\overline{\Omega})$ and $u\in C(\overline{Q_T})\cap
C^{1+\nu/2, 2+\nu}(\overline{Q_{\tau,T}})$ the solution of problem
\eq{part3.1.1} found above. Consider $w_0\in C(\overline{\Omega})$ and $w\in C(\overline{Q_T})\cap
C^{1+\nu/2, 2+\nu}(\overline{Q_{\tau,T}})$ a function satisfying:
\begin{equation}\label{part3.7.1}
\left\{\begin{array}{l}
\di\frac{\partial w}{\partial t}-\Delta w\leq a|\nabla w|^p\ (\ \geq a|\nabla
w|^p\ )\ \ \mbox{in} \ \ Q_T,\\
\di\frac{\partial w}{\partial\nu}\leq 0\ (\ \geq 0\ )\ \ \mbox{on}\ \
\Gamma_T,\\
w(0,\cdot)=w_0\leq\mu_0\ (\mbox{resp.}\ w_0\geq\mu_0\ )\ \ \mbox{in}\ \ \Omega.\end{array}\right.
\end{equation}
Then:
$$w\leq u \ (\mbox{resp.}\ w\geq u\ )\ \ \mbox{in}\ \ Q_T.$$
\end{lemma}
\proof
An analogous result for the whole space $\rr^N$ can be found in \cite{GGK}
[Lemma 7] and our proof follows the same arguments.\\
We suppose first that $\Omega$ is a bounded and convex open set which
contains the origin and $w_0\leq \mu_0$. \\
Consider two real numbers $\eps>0$ and  $A>0$, and denote by $z$ the function:
\begin{equation}\label{part3.part2.4.11}
z(t,x)=w(t,x)-u(t,x)-At^{q}-\eps(1+|x|^2)^{\demi},
\end{equation}
where
 $q=\min\{1,\frac{1}{p}\}.$
Then:
$$z\in C(\overline{Q_T})\cap C^{1,2}((0,T]\times\overline{\Omega})\ \
\mbox{and}\ \ z(0,x)\leq 0 \ ~\mbox{for all} \ ~ x\in\Omega.$$
Thanks to Lemma \pr{part3.l.3} and to hypothesis \eq{part3.7.1} we have:
\begin{equation}\label{part3.part2.4.12}
\frac{\partial z}{\partial \nu}(t,x)=\frac{\partial w}{\partial \nu}(t,x)-\frac{\partial u}{\partial \nu}(t,x)-\frac{\eps\cdot
  x\cdot\nu}{\sqrt{1+|x|^2}}<0 \ ~\mbox{on} \ ~\Gamma_T.
\end{equation}
We claim that
\begin{equation}\label{part3.part2.4.13}
z(t,x)\leq 0 \ ~\mbox{for all} \ ~(t,x)\in Q_T,
\end{equation}
Indeed, if $z$ is positive anywhere in $Q_T$
then $z$ has a positive maximum in $(t_0,x_0)\in (0,T]\times\Omega$ since, if
$(t_0,x_0)\in (0,T]\times\partial\Omega$, we have:
$$\frac{\partial z}{\partial\nu}(t_0,x_0)=\lim\limits_{\lambda\nearrow
0}\frac{z(t_0,x_0+\lambda\nu)-z(t_0,x_0)}{\lambda}\geq 0,$$
which contradicts relation \eq{part3.part2.4.12}.The rest of the proof is standard and follows the same ideas as the proof
of Lemma 7 in \cite{GGK}. So, it will be omitted.

In the general case,when $\Omega$ does not
contain the origin it is enough to
translate the problem on a domain which contains the origin, for example
$\Omega_{x_0}=\Omega-x_0$ where $x_0\in\Omega$.
\endproof
\begin{rem}\label{part1.r.10}
The result of Lemma \pr{part3.l.6} is valid for all $a\in\rr, a\neq 0$
an so, for any solution $u\in C(\overline{Q_T})\cap
C^{1+\delta/2,2+\delta}(\overline{Q_{\tau,T}})$ of problem \eq{part3.1.1}.
\end{rem}

\section{Proof of Theorem \pr{part3.t.2}}
Let $\beta $ be a positive number satisfying:
\begin{equation}\label{part3.4.1}
\beta>\left(\frac{2p+1-N}{N}\right)_{+}
\end{equation}
and set:
\begin{equation}\label{part3.4.2}
\gamma=N(\beta+1)\ \ \ ~\mbox{and}\ \ \ ~\eta=N(\beta+1)-p.
\end{equation}
Since $\beta$ satisfies \eq{part3.4.1} we have:
\begin{equation}\label{part3.4.3}
\eta>p+1
\end{equation}
with this notations, we can state the following proposition which is the key argument in the
proof of Theorem \pr{part3.t.2}.
\begin{prop}\label{part3.p.3}
Let $\Omega\subset\rr^N$ be a bounded and convex domain with smooth
boundary and $\mu_0\in C(\overline{\Omega})$. Let denote by  $u$ the
solution of problem \eq{part3.1.1} whose existence was proved in Theorem \pr{part3.t.1}. Then:
\begin{itemize}
\item[$i)$]  The application
$t\mapsto(1+t)(M(t)-m(t))^{\gamma}$ belongs to $L^1(0,+\infty)$,
\item[$ii)$] Denoting by $y$ the function defined on $[0,+\infty)$ by:
$$y(t)=\int\limits_{t}^{\infty}(s-t)(M(s)-m(s))^{\gamma}\, ds,$$
then $y\in W^{2,\infty}((0,+\infty))$ and satisfies the following
differential inequality:
\begin{equation}\label{part3.4.4}
y'(t)+C y(t)^{\alpha}\leq 0,\ \ \forall \ t\in[0,+\infty),
\end{equation}
where the positive constant depends only on  $p,\beta,N, \Omega, (M(0)-m(0))$ and
\begin{equation}\label{part3.4.5}
\alpha=\frac{1+\eta}{2+\eta-p}.
\end{equation}
\end{itemize}
\end{prop}
\proof
The proof follows the same ideas as those of  Lemma 3 in
\cite{BLS}  and Lemma 12 in \cite{BLSS}. Setting:
$$T^{*}=\inf\{t>0; |\nabla u(t)|\equiv  0\}$$
then $T^*$ can be also defined by:
\begin{equation}\label{part3.4.16}
T^*=\inf\{t>0; M(t)=m(t)\}=\inf\{t>0; y(t)=0\}.
\end{equation}
First, if $\mu_0\equiv c$ then $u\equiv
c$, which implies $T^*=0$ and Proposition \pr{part3.p.3} is achieved. \\
We suppose that $\mu_0$ is not constant, consequently
$T^*\in(0,+\infty]$. Consider $T\in(0,T^*)$ and $t\in[0,T)$. Integrating
the first equation of problem \eq{part3.1.1} on $(t,T)\times\Omega$, and using
relations \eq{part3.1.2}, \eq{part3.1.3}, \eq{part3.1.4} and
\eq{part3.1.5} we get:
$$|a|\int\limits_{t}^{T}\int_{\Omega}|\nabla u(s,x)|^p\,dx\,ds\leq |\Omega|
(M(t)-m(t))$$
Recalling \eq{part3.4.2}  and \eq{part3.4.3} we have:
$\gamma=N(\beta+1)=\eta+p$
and we deduce that:
$$\|\nabla u(s)\|^{\gamma}_{\gamma}\leq\|\nabla
u(s)\|^{\eta}_{\infty}\cdot\||\nabla u(s)|^p\|_{1}.$$ Combining
these two last inequalities we get:
\begin{equation}\label{part3.4.6}
\int\limits_{t}^{T}\|\nabla u(s)\|^{-\eta}_{\infty}\cdot\|\nabla
u(s)\|^{\gamma}_{\gamma} ds\leq\frac{|\Omega|}{|a|}( M(t)-m(t))
\end{equation}

We distinguish between the two cases below:\\
$(i)$ The case $p\neq 1$.
Thanks to relation \eq{part3.1.6}, for all $s\in(t,T)$ we have:
\begin{equation}\label{part3.4.7}
\|\nabla u(s)\|_{\infty}^{-\eta}\geq C_1(s-t)^{\eta/p}(M(t)-m(t))^{-\eta/p},
\end{equation}
where  $C_1$ is a positive constant which depends only on $p>0$.\\
Applying Lemma \pr{part3.l.2} we have:
\begin{equation}\label{part3.4.8}
M(s)-m(s)\leq C_2\|\nabla u(s)\|_{\gamma}
\end{equation}
And combining \eq{part3.4.6},\eq{part3.4.7} and \eq{part3.4.8} we get:
\begin{equation}\label{part3.4.9}
\int\limits_{t}^{T}(s-t)^{\eta/p}(M(s)-m(s))^{\gamma}\, ds\leq
C_3(M(t)-m(t))^{\gamma/p}\ \mbox{for all}\ t\in[0,T],
\end{equation}
where the constant $C_3$ depends only on $N, a, p, \beta,\eta$ and
$\Omega$.
\\
 We can pass to the limit in  \eq{part3.4.9}, as $T\nearrow T^{*}$,  and we
 obtain:
\begin{equation}\label{part3.4.10}
\int\limits_{t}^{T^*}(s-t)^{\eta/p}(M(s)-m(s))^{\gamma}\,ds\leq
C_3(M(t)-m(t))^{\gamma/p}\ \mbox{for all}\ t\in[0,T^*)
\end{equation}
Let fix $\delta \in(0,T^*)$. Using
\eq{part3.1.4},\eq{part3.1.5},\eq{part3.4.3} and \eq{part3.4.10} we get:
$$\begin{array}{l}
\di\int\limits_0^{\infty}(1+s)(M(s)-m(s))^{\gamma}\,ds=\\
\ \ \ \ \ \ \ =\di\int\limits_0^{\delta}(1+s)(M(s)-m(s))^{\gamma}\,
  ds+\di\int\limits_{\delta}^{T^*}(1+s)(M(s)-m(s))^{\gamma}\,ds\ \ \ \ \\
\ \ \ \ ~\leq (1+\delta)[\delta(M(0)-m(0))^{\gamma}+\delta
^{-\eta/p}\di\int\limits_{\delta}^{T^*}s^{\eta/p}(M(s)-m(s))^{\gamma}\,ds]\\
\ \ \ \ ~\leq (1+\delta)[\delta(M(0)-m(0))^{\gamma}+\delta
^{-\eta/p}\di\int\limits_{0}^{T^*} s^{\eta/p}(M(s)-m(s))^{\gamma}\,ds].
\end{array}$$
Once again, using \eq{part3.1.4}, \eq{part3.1.5} and \eq{part3.4.10} (which
is, in particular, valid for $t=0$), this last integral is finite. Consequently:
\begin{equation}\label{part3.4.11}
t\mapsto(1+t)(M(t)-m(t))^{\gamma}\in L^1((0,+\infty)).
\end{equation}
And we deduce that the function $y$ is well defined on
$[0,+\infty)$ and belongs to
$W^{2,\infty}((0,+\infty))$. Indeed, for $t>0$, we have:
\begin{equation}\label{part3.4.11.1}
y'(t)=-\int\limits_t^{+\infty}(M(s)-m(s))^{\gamma}\,ds
\end{equation}
and
\begin{equation}\label{part3.4.11.2}
y^{''}(t)=(M(t)-m(t))^{\gamma}.
\end{equation}
Using  H\"older inequality and \eq{part3.4.3} we deduce that the function $y$ verifies:
$$\begin{array}{l}
\di(y(t))^{\eta/p}=(\int\limits_t^{\infty}(s-t)(M(s)-m(s))^{\gamma}\,ds)^{\eta/p}\\
\di\ \ \ \ \ \ \ \ \ \ \ \ \ \leq
\Big[\int\limits_{t}^{T^*}(s-t)^{\eta/p}(M(s)-m(s))^{\gamma}\,ds\Big]\Big[\int\limits_t^{T^*}(M(s)-m(s))^{\gamma}\,ds\Big]^{\eta/p-1}\end{array}$$
Combining this last inequality with \eq{part3.4.9} we get:
\begin{equation}\label{part3.4.12}
y(t)^{\eta/p}\leq C\cdot y''(t)^{1/p}(-y'(t))^{\eta/p-1},
\end{equation}
which yields:
\begin{equation}\label{part3.4.13}
y(t)^{\eta}\leq C\cdot y''(t)\cdot (-y'(t))^{\eta-p},\ \ \forall\ t\in[0,T^*),
\end{equation}
Taking into account the fact that $y'(t)\leq 0$, we can multiply \eq{part3.4.13} by $(-y'(t))$
and integrate over $(t,T^*)$. We get:
$$y(t)^{1+\eta}\leq C\cdot (-y'(t))^{2+\eta-p},\ \forall\  t\in[0,T^*),$$
and thanks to the definition of $T^*$ it follows that:
$$y'(t)+\frac{1}{C}y(t)^{\frac{1+\eta}{2+\eta-p}}\leq 0, \ \ \forall \ t\in[0,+\infty).$$
Hence \eq{part3.4.4} holds for $p\neq 1$.\\
$(ii)$ The case $p=1$. Instead of \eq{part3.1.6} we can use this time
\eq{part3.1.5.1}. Thus, for all $s\in (t,T)$, we have:
\begin{equation}\label{part3.4.14}
\|\nabla u(s)\|_{\infty}^{-\eta}\geq
C_4(s-t)^{\eta/2}(M(t)-m(t))^{-\eta}.
\end{equation}
Combining relations \eq{part3.4.6}, \eq{part3.4.14} and \eq{part3.4.8} we get:
\begin{equation}\label{part3.4.15}
\int\limits_{t}^{T}(s-t)^{\eta/2}(M(s)-m(s))^{\gamma}\, ds\leq
C_5 (M(t)-m(t))^{\gamma}\ \mbox{for all}\ t\in[0,T]
\end{equation}
where $C_5$ is a positive constant which depends only on $N, a, p, \beta,\eta$ and
$\Omega$.\\
Thanks to \eq{part3.4.3} we have $\eta>2$. As previously, we fix $\delta
\in(0,T^*)$, then, using \eq{part3.1.4}, \eq{part3.1.5}, \eq{part3.4.3} and
\eq{part3.4.15} we get:
$$\begin{array}{l}
\di\int\limits_0^{\infty}(1+s)(M(s)-m(s))^{\gamma}\,ds=\\
\ \ \ \ \ \ =\di\int\limits_0^{\delta}(1+s)(M(s)-m(s))^{\gamma}\,
  ds+\di\int\limits_{\delta}^{T^*}(1+s)(M(s)-m(s))^{\gamma}\,ds\\
\ \ \ \ ~\leq (1+\delta)[\delta(M(0)-m(0))^{\gamma}+\delta
^{-\eta/2}\di\int\limits_{0}^{T^*} s^{\eta/2}(M(s)-m(s))^{\gamma}\,ds].
\end{array}$$
>From \eq{part3.1.4}, \eq{part3.1.5} and \eq{part3.4.15}, this last
integral is finite. Consequently relation \eq{part3.4.11} is valid for $p=1$,
too.
As in the first case we deduce that the function $y$ is well defined on
$[0,+\infty)$ and belongs to
$W^{2,\infty}((0,+\infty))$, the first and the second derivatives being
given by \eq{part3.4.11.1} and
\eq{part3.4.11.2}. Since $\eta>2$, using H\"older inequality we get this time:
$$\begin{array}{l}
\di(y(t))^{\eta/2}=(\int\limits_t^{\infty}(s-t)(M(s)-m(s))^{\gamma}\,ds)^{\eta/2}\\
\di\ \ \ \ \ \ \ \ \ \ \ \ \ \ \ \ \leq
\Big[\int\limits_{t}^{T^*}(s-t)^{\eta/2}(M(s)-m(s))^{\gamma}\,ds\Big]\Big[\int\limits_t^{T^*}(M(s)-m(s))^{\gamma}\,ds\Big]^{\eta/2-1}\end{array}$$
Taking into account \eq{part3.4.15}, \eq{part3.4.11.1} and
\eq{part3.4.11.2} we deduce:
\begin{equation}\label{part3.4.12.1}
y(t)^{\eta/2}\leq C\cdot y''(t)(-y'(t))^{\eta/2-1}
\end{equation}
and by the same arguments as previously we get:
$$y'(t)+\frac{1}{C}y(t)\leq 0 \ \mbox{for all} \ t\in[0,+\infty)$$
which ends the proof of Proposition \pr{part3.p.3}, as $p=1$.
\endproof
{\bf Proof of Theorem \pr{part3.t.2}:}
Let:
$$y(t)=\int\limits_t^{\infty}(s-t)(M(s)-m(s))^{\gamma}\,ds,\ ~t\in[0,\infty),$$
be the function defined in Proposition \pr{part3.p.3}. We have obtained that $y\in
W^{2,+\infty}((0,+\infty))$ and there is a positive constant $C$
depending only on $p,\beta,N,\Omega$ and $(M(0)-m(0))$ such that $y$
satisfies the differential inequality \eq{part3.4.4}:
$$y'(t)+C y(t)^{\alpha}\leq 0,\ \mbox{for all}\ t\in[0,+\infty),$$
with  $\alpha$ given by \eq{part3.4.5} and
$y(0)=\int\limits_0^{\infty}s(M(s)-m(s))^{\gamma}\,ds \geq 0.$

On the one hand if $p\in(0,1)$ then $\alpha\in(0,1)$ and thanks to
\eq{part3.4.4} and \eq{part3.4.16} we get that $T^*<\infty$ and:
$$y(t)\equiv 0\ \mbox{for}\ t>T^*$$
Consequently, for $p\in(0,1)$,  the extinction of the gradient in finite time of the
solution to problem \eq{part3.1.1}  occurs.

On the other hand, if $p\geq 1$ then $\alpha\geq 1$ and thanks to
\eq{part3.4.4} and from the fact that $y$ is a positive function, we deduce that:
\begin{equation}\label{part3.5.100}
\frac{y'(t)}{y^{\alpha}(t)}\leq  -\frac{1}{C}
 \ \ \mbox{for all}\ \ t\in(0,+\infty).
\end{equation}
We distinguish between the two cases below:\\
$(i)$ The case $p>1$. We have $\alpha>1$ and integrating \eq{part3.5.100}, over
 $(0,t)$, $t>0$, we obtain:
$$y(t)^{1-\alpha}\geq y(0)^{1-\alpha}+\frac{(\alpha-1)t}{C},$$
or else:
$$y(t)\leq\left(\frac{1}{y(0)^{1-\alpha}+\frac{(\alpha-1)t}{C}}\right)^{1/(\alpha-1)}.$$
$(ii)$ The case $p=1$. We have $\alpha=1$ and integrating \eq{part3.5.100}
over $(0,t)$, $t>0$ we get this time:
$$\log y(t)\leq \log y(0)-\frac{t}{C},$$
or else:
$$y(t)\leq y(0)e^{-\frac{t}{C}}.$$
Thus, we have obtained the decreasing rate for the function $y$, as
$t\to\infty$. \\
We claim that:
\begin{equation}\label{part3.4.20}
\lim\limits_{t\to\infty}M(t)-m(t)=0,
\end{equation}
which implies that, for $p\geq 1$, the solution $u$ of problem \eq{part3.1.1}
converges uniformly in $\overline{\Omega}$ to a constant, as
$t\to\infty$. \\
To prove \eq{part3.4.20} we recall that,  the
function defined by:
\begin{center}
$g(t)=(1+t)(M(t)-m(t))^{\gamma}$ belongs to $L^1(0,+\infty)$.
\end{center}
Since $t\mapsto(M(t)-m(t))$ is a positive and decreasing function on
$[0,+\infty)$, there exists a positive constant $c$ such that:
$$c=\lim\limits_{t\to\infty}(M(t)-m(t))^{\gamma}.$$
Then, $c(1+t)\leq g(t),$ and it follows that the function $t\to
c(1+t)$ belongs to $L^1(0,+\infty)$, which is possible only if
$c=0$. So assertion \eq{part3.4.20} holds. Now, we want to find a
decreasing rate for the application $t\mapsto (M(t)-m(t))$. Denote
by $f$ the decreasing rate of the function  $y$:
\begin{equation}\label{part3.5.101}
f(t)=\left\{\begin{array}{l}
y(0)e^{-\frac{t}{C}}\ \ \mbox{if} \ \ p=1,\\
\left(\frac{1}{y(0)^{(1-\alpha)}+\frac{(\alpha-1)t}{C}}\right)^{1/(\alpha-1)}\
\ \mbox{if} \ \ p>1.
\end{array}\right.
\end{equation}
Then, for all $t>0$ we have:
$$\int\limits_{t/2}^t (s-t/2)(M(s)-m(s))^{\gamma}\,ds\leq y(t/2)\leq f(t/2).$$
Since $s\mapsto (M(s)-m(s))$ is a decreasing function, we deduce that
$$(M(t)-m(t))^{\gamma}\int\limits_{t/2}^t (s-t/2)\, ds\leq f(t/2),\ \
\forall\ t>0,$$
which implies:
\begin{equation}\label{part3.5.102}
M(t)-m(t)\leq \left(\frac{8f(t/2)}{t^2}\right)^{\frac{1}{\gamma}},\ \
\forall\ t>0,
\end{equation}
where  $f$ is given by \eq{part3.5.101}. This ends the proof of Theorem \pr{part3.t.2}.
\endproof

\end{document}